\documentclass[11pt]{article}
\usepackage[english]{babel}
\usepackage{amsmath,amsthm,amssymb}
\usepackage{mathrsfs}
\usepackage{bbm}
\usepackage{amsfonts}
\usepackage[margin=0.8in]{geometry}
\usepackage{graphicx}
\usepackage{caption}
\usepackage{subcaption}
\usepackage{wasysym}
\usepackage{enumerate}
\usepackage{algorithm2e}
\usepackage{tabularx}
\usepackage{color}
\usepackage{wrapfig}
\usepackage{todonotes}

\newtheorem{thm}{Theorem}[section]
\newtheorem{ass}[thm]{Assumption}
\newtheorem{cor}[thm]{Corollary}
\newtheorem{lem}[thm]{Lemma}
\newtheorem{prop}[thm]{Proposition}
\newtheorem*{hyp*}{Hypothesis}
\theoremstyle{definition}
\newtheorem{defn}[thm]{Definition}

\theoremstyle{remark}
\newtheorem{rem}[thm]{Remark}

\numberwithin{equation}{section}
\newcommand{\R}{\mathbb R}
\newcommand{\eps}{\varepsilon}

\newcommand{\bbF}{\mathbb F}

\newcommand{\bbG}{\mathbb G}

\newcommand{\bbN}{\mathbb N}

\newcommand{\mcA}{\mathcal{A}}

\newcommand{\mcB}{\mathcal{B}}

\newcommand{\mcE}{\mathcal E}
\newcommand{\mcG}{\mathcal{G}}

\newcommand{\mcF}{\mathcal F}

\newcommand{\mcL}{\mathcal L}
\newcommand{\mcT}{\mathcal T}
\newcommand{\mcP}{\mathcal P}

\newcommand{\mcH}{\mathcal H}
\newcommand{\mcK}{\mathcal K}
\newcommand{\mcU}{\mathcal U}

\newcommand{\mcS}{\mathcal S}

\newcommand{\mcV}{\mathcal V}
\newcommand{\mcY}{\mathcal Y}
\newcommand{\mcZ}{\mathcal Z}

\newcommand{\E}{\mathbb{E}}
\newcommand{\Prob}{\mathbb{P}}

\newcommand{\esssup}{\mathop{\rm{ess}\,\sup}}
\newcommand{\essinf}{\mathop{\rm{ess}\,\inf}}

\newcommand{\ett}{\mathbbm{1}}

\newcommand{\cadlag}{c\`adl\`ag~}

\newcommand{\Pred}{\mathcal{P}}

\newcommand{\trace}{{\rm Tr}}

\newcommand{\bigS}{\mathfrak S}

\newcommand{\infOP}{\mathcal M}
\newcommand{\supOP}{\mathcal W}

\newcommand{\ie}{\textit{i.e.\ }}
\newcommand{\eg}{\textit{e.g.\ }}

\begin{document}

\title{Probabilistic Representation for Viscosity Solutions to Double-Obstacle Quasi-Variational Inequalities}

\author{Magnus Perninge\footnote{M.\ Perninge is with the Department of Physics and Electrical Engineering, Linnaeus University, V\"axj\"o,
Sweden. e-mail: magnus.perninge@lnu.se.}} %
\maketitle
% ----------------------------------------------------------------
\begin{abstract}
We prove the existence and uniqueness of viscosity solutions to quasi-variational inequalities (QVIs) with both upper and lower obstacles. In contrast to most previous works, we allow all involved coefficients to depend on the state variable and do not assume any type of monotonicity. It is well known that double obstacle QVIs are related to zero-sum games of impulse control, and our existence result is derived by considering a sequence of such games. Full generality is obtained by allowing one player in the game to randomize their control. A by-product of our result is that the corresponding zero-sum game has a value, which is a direct consequence of viscosity comparison.

Utilizing recent results for backward stochastic differential equations (BSDEs), we find that the unique viscosity solution to our QVI is related to optimal stopping of BSDEs with constrained jumps and, in particular, to the corresponding non-linear Snell envelope. This gives a new probabilistic representation for double obstacle QVIs.

It should be noted that we consider the min-max version (or equivalently the max-min version); however, the conditions under which solutions to the min-max and max-min versions coincide remain unknown and is a topic left for future work.
\end{abstract}

% ----------------------------------------------------------------
\section{Introduction}
We establish existence and uniqueness of viscosity solutions to double obstacle quasi-variational inequalities (QVIs) of the type
\begin{align}\label{ekv:var-ineq}
\begin{cases}
 \min\{v(t,x)-\supOP v(t,x), \max\{v(t,x)-\infOP v(t,x),-\frac{\partial}{\partial t}v(t,x)-\mcL v(t,x)\\
  \quad-f(t,x,v(t,x),\sigma^\top(t,x)\nabla_x v(t,x))\}\}=0,\quad\forall (t,x)\in[0,T)\times \R^d \\
  v(T,x)=\psi(x),
\end{cases}
\end{align}
where $\supOP v(t,x):=\sup_{b\in U}\{v(t,x+\xi(t,x,b))-\ell(t,x,b)\}$ and $\infOP v(t,x):=\inf_{e\in E}\{v(t,x+\gamma(t,x,e))+\chi(t,x,e)\}$ for $t\in [0,T)$ and
\begin{align}
  \mcL:=\sum_{j=1}^d a_j(t,x)\frac{\partial}{\partial x_j}+\frac{1}{2}\sum_{i,j=1}^d (\sigma\sigma^\top(t,x))_{i,j}\frac{\partial^2}{\partial x_i\partial x_j}
\end{align}
is the infinitesimal generator related to the SDE
\begin{align*}
\check X^{t,x}_s=x + \int_t^s a(r,\check X^{t,x}_r)dr + \int_t^s \sigma(r,\check X^{t,x}_r)dW_r,\quad\forall s\in[t,T]
\end{align*}
driven by a $d$-dimensional Brownian motion, $W$.

Existence of a viscosity solution to \eqref{ekv:var-ineq} is derived by considering a zero-sum game of impulse control where one player (the player that has act-first advantage) implements a purely randomized impulse control, avoiding the possibility for preemption, while the opponent chooses a standard impulse control. For each $(t,x)\in [0,T]\times\R^d$ and any impulse control $u:=(\tau_j,\beta_j)_{j=1}^N\in\mcU_t$ (the subscript $t$ is used to indicate that $\tau_1\geq t$), we let the controlled process $X^{t,x;u}$ solve the SDE with jumps and impulses,
\begin{align}\nonumber
X^{t,x;u}_s&=x+\int_t^s a(r,X^{t,x;u}_r)dr+\int_t^s\sigma(r,X^{t,x;u}_r)dW_r
\\
&\quad + \int_t^s\!\!\!\int_E\gamma(r,X^{t,x;u}_{r-},e)\mu(dr,de)+\sum_{j=1}^{N}\ett_{[\tau_j\leq s]}\xi(\tau_j,X_{\tau_j}^{t,x;[u]_{j-1}},\beta_j),\quad\forall s\in [0,T],\label{ekv:fwd-sde-Markov}
\end{align}
where $\mu$ is a Poisson random measure on $[0,T]\times E$, with a compensator $\lambda$ that has full topological support. The random measure, $\mu$, models randomized impulses and the randomization is controlled by a change of measure that effectively changes its jump intensity. The accumulated impulse cost for the second player gives rise to a process $\Xi$, defined as
\begin{align}\label{ekv:Xi-def}
\Xi^{t,x;u}_s:=\sum_{j=1}^N\ett_{[\tau_j<s]}\ell(\tau_j,X^{t,x;[u]_{j-1}}_{\tau_{j}},\beta_j).
\end{align}
Given an admissible density $\nu\in\mcV_t$, \ie a predictably-measurable bounded map $\nu:[t,T]\times\Omega\times E\to [0,\infty)$, we let $\mu^\nu(ds,de):=\mu(ds,de)-\nu_s(e)\lambda(de)ds$ and consider the BSDE
\begin{align}\nonumber
P^{t,x;u,\nu}_s&=\psi(X^{t,x;u}_T)+\int_s^T f(r,X^{t,x;u}_r,P^{t,x;u,\nu}_r,Q^{t,x;u,\nu}_r)dr-\int_s^T Q^{t,x;u,\nu}_r dW_r
\\
&\quad-\int_{s}^T\!\!\!\int_E R^{t,x;u,\nu}_{r}(e)\mu^\nu(dr,de) + \int_{s}^T\!\!\!\int_E\chi(r,X^{t,x;u}_{r},e)\nu_r(e)\lambda(de)dr - \Xi^{t,x;u}_{T}+\Xi^{t,x;u}_s. \label{ekv:non-ref-bsde-simp-alt}
\end{align}
Using the unique solution to \eqref{ekv:non-ref-bsde-simp-alt}, we define the lower (resp. upper) value function as
\begin{align*}
  V^-(t,x):=\sup_{u\in\mcU^{t,x}_t}\inf_{\nu\in\mcV^{t,x}}P^{t,x;u,\nu}_t\quad(\text{resp. } V^+(t,x):=\inf_{\nu\in\mcV^{t,x}}\sup_{u\in\mcU^{t,x}_t}P^{t,x;u,\nu}_t),
\end{align*}
where $\mcU^{t,x}_t$ (resp.~$\mcV^{t,x}$) is the set of impulse controls (resp.~admissible densities) measurable with respect to a filtration generated by the process $X^{t,x;u}$. We demonstrate that the partial differential equation (PDE) in \eqref{ekv:var-ineq} is closely related to the zero-sum game discussed above. By applying two different types of truncations to the game, we extract a supersolution $\underline v$ and a subsolution $\bar v$ to \eqref{ekv:var-ineq}, respectively, that trivially satisfy $\underline v(t,x)\leq V^-(t,x)$ and $V^+(t,x)\leq \bar v(t,x)$. Value of the game then follows by a comparison result for viscosity solutions to \eqref{ekv:var-ineq} that requires an additional no-free-loop type of condition.

From the seminal work presented in \cite{Kharroubi2010} (see also \cite{Bouchard09}), it is well known that control randomization of the type described above is connected to BSDEs with a constraint on the jump process. We further extend this machinery by deriving an alternative representation for the unique viscosity solution to \eqref{ekv:var-ineq}. In particular, we relate $v$ to the non-linear Snell envelope
\begin{align}\label{ekv:Snell-env}
Y^{t,x}_s:=\esssup_{\tau\in\mcT_s}Y^{t,x;\tau}_s,
\end{align}
where for each stopping time $\tau\geq t$, the process $Y^{t,x;\tau}$ is the first component in the quadruple of processes $(Y^{t,x;\tau},Z^{t,x;\tau},V^{t,x},K^{-,t,x;\tau})$ that is the unique maximal solution (in the sense that $Y^{t,x;\tau}\geq \tilde Y^{t,x;\tau}$ for any other solution $\tilde Y^{t,x;\tau}$) to
\begin{align}\label{ekv:bsde-c-jmp}
  \begin{cases}
    Y^{t,x;\tau}_s=\supOP Y(\tau,X^{t,x}_\tau)+\int_s^\tau f(r,X^{t,x}_r,Y^{t,x;\tau}_r,Z^{t,x;\tau}_r)dr-\int_s^\tau Z^{t,x;\tau}_r dW_r-\int_{s}^\tau\!\!\!\int_E V^{t,x;\tau}_{r}(e)\mu(dr,de)
    \\
    \quad -(K^{-,t,x;\tau}_\tau-K^{-,t,x;\tau}_s),\quad \forall s\in [t,\tau],
    \\
    V^{t,x;\tau}_s(e)\geq - \chi(s,X^{t,x}_{s-},e),\quad d\Prob\otimes ds\otimes \lambda(de)-a.e.
  \end{cases}
\end{align}
Here, $Y(t,x):=Y^{t,x}_t$ and $X^{t,x}$ is the unique solution to the SDE
\begin{align}\label{ekv:fwd-sde}
X^{t,x}_s&=x+\int_t^s a(r,X^{t,x}_r)ds+\int_0^t\sigma(r,X^{t,x}_r)dW_r + \int_t^s\!\!\!\int_E\gamma(r,X^{t,x}_{r-},e)\mu(dr,de),\quad\forall s\in[0,T].
\end{align}

%Our approach is based on a combination of truncation and penalization and we introduce the double sequence of penalized BSDEs with jumps
%\begin{align}\label{ekv:rbsde-pen-trunk}
%  \begin{cases}
%    Y^{t,x,k,n}_s=\psi(X^{t,x}_T)+\int_s^T f(r,X^{t,x}_r,Y^{t,x,k,n}_r,Z^{t,x,k,n}_r)dr-n\int_t^T\!\!\int_E(V^{t,x,k,n}_r(e)+\chi(r,X^{t,x}_{r},e))^-\lambda(de)dr
%    \\
%    \quad -\int_s^T Z^{t,x,k,n}_r dW_r-\int_{s}^T\!\!\!\int_E V^{t,x,k,n}_{r}(e)\mu(dr,de)+K^{t,x,k,n}_T-K^{t,x,k,n}_s,\quad \forall s\in [t,T],
%    \\
%    Y^{t,x,k,n}_s\geq \supOP Y^{k-1,n}(s,X^{t,x}_s),\, \forall s\in [t,T] \quad\text{and}\quad\int_t^T\!\! \big(Y^{t,x,k,n}_s-\supOP Y^{k-1,n}(s,X^{t,x}_s)\big)dK^{t,x}_s
%  \end{cases}
%\end{align}
%where $Y(t,x):=Y^{t,x}_t$ and set $K^{-,n}_t:=n\int_0^t\!\!\int_E(V^{n}_s(e)+\chi(s,Y^{n}_{s},Z^{n}_s,e))^-\lambda(de)ds$.
%
%Moreover, we introduce the

\textbf{Related literature} When the lower barrier is absent (\ie when $\ell\equiv\infty$), the PDE described in \eqref{ekv:var-ineq} reduces to a standard (single-barrier) QVI. It is well known that, under suitable conditions on the involved parameters, value functions of impulse control problems are solutions (in viscosity sense) to standard QVIs when the driving noise process is a Brownian motion (see the seminal work in~\cite{BensLionsImpulse}) and to so-called quasi-integrovariational inequalities when the driving noise is a general L\'evy process~\cite{OksenSulemBok}. Employing a contraction argument, QVIs featuring general non-local drivers were investigated in \cite{qvi-rbsde}. Specifically, \cite{qvi-rbsde} establishes a Feynman-Kac representation of such QVIs by linking their solutions to systems of reflected BSDEs (RBSDEs).

In the case when the impulse costs $\ell$ and $\chi$ are non-increasing in time and independent of the state variable, $x$, a zero-sum game of impulse control was resolved in \cite{Cosso13} by proving that the corresponding value function is indeed the unique viscosity solution to a double obstacle QVI. The results presented in \cite{Cosso13} extended previous results from \cite{TangHou07}, where one player was restricted to play a switching control. Since then, zero-sum games where both players implement switching controls have been considered in a more general framework in \cite{BollanSWG1} and \cite{DjehicheSWG2}.

An alternative Feynman-Kac representation for solutions to standard QVIs was proposed in \cite{Kharroubi2010} (see also \cite{Bouchard09}), where the solution to a QVI is related to the minimal solution of a BSDE with constrained jumps. Following the seminal work in \cite{Kharroubi2010}, BSDEs with positive jumps were related to fully non-linear Hamilton-Jacobi-Bellman integro-partial differential equations (HJB-IPDEs) through a Feynman-Kac representation in \cite{KharroubiPham15}, while a correspondence between RBSDEs with positive jumps and fully non-linear variational inequalities was established in \cite{Choukroun15}.

The ensemble of approaches to find probabilistic representations of PDEs or to solve stochastic optimal control problems that utilize BSDEs with constrained jumps is commonly referred to as \emph{control randomization}. Recently, \cite{imp-stop-game} further extended the scope of control randomization by establishing a link between the solution to a reflected BSDE with constrained jumps and optimal stopping of BSDEs with constrained jumps. The results in \cite{imp-stop-game} are applicable to the optimal stopping problem described in \eqref{ekv:Snell-env}-\eqref{ekv:fwd-sde}, provided that $\supOP Y(\tau,X^{t,x}_\tau)$ is replaced by a regular barrier $h(\tau,X^{t,x}_\tau)$. The corresponding Feynman-Kac representation was derived in \cite{qvi-stop}, which relates $Y^{t,x}_s$ to the unique viscosity solution to \eqref{ekv:var-ineq}, with $h(t,x)$ in place of $\supOP v(t,x)$. Notable is that the approach taken in the present work relies heavily on the developments in \cite{imp-stop-game} and \cite{qvi-stop}.\\

\textbf{Outline} In the next section, we set the notations and state the assumptions that hold throughout. In addition, we give some preliminary results that are repeatedly referred in the article. In Section~\ref{sec:game}, we formulate the zero-sum game and state our main result (Theorem~\ref{thm:game-value}). In the following section, we turn to the upper and lower bounds $\underline v$ and $\bar v$ and show that $\underline v$ (resp.~$\bar v$) is a supersolution (resp.~subsolution) to \eqref{ekv:var-ineq}. Together with the comparison principle derived in Appendix~\ref{app:uniq} (effectively stating that a supersolution to \eqref{ekv:var-ineq} always dominates a subsolution) this proves existence of a unique viscosity solution to \eqref{ekv:var-ineq}, along with a value for the game. Finally, in Section~\ref{sec:prob-rep} we turn to the probabilistic representation which, in light of our earlier findings, is derived in a straightforward manner.

\section{Preliminaries\label{sec:prel}}
\subsection{Notations}
We let $(\Omega,\mcF,\Prob)$ be a complete probability space on which lives a $d$-dimensional Brownian motion $W$ and an independent Poisson random measure $\mu$ defined on $[0,T]\times E$ with intensity $\lambda(de)$. Here, it is assumed that $E$ is a compact subset of $\R^d$ endowed with its Borel $\sigma$-field $\mcB(E)$. We denote by $\bbF:=(\mcF_t)_{0\leq t\leq T}$ the augmented natural filtration generated by $W$ and $\mu$ and for $t\in[0,T]$ we let $\bbF^t:=(\mcF^t_s)_{t\leq s\leq T}$ denote the augmented natural filtration generated by $(W_s-W_t:t\leq s\leq T)$ and $\mu(\cdot\cap[t,T],\cdot)$.\\

\noindent Throughout, we will use the following notations:
\begin{itemize}
\item We let $\Pi^g$ denote the set of all functions $\varphi:[0,T]\times\R^d\to\R$ that are of polynomial growth in $x$, \ie there are constants $C,\rho>0$ such that $|\varphi(t,x)|\leq C(1+|x|^\rho)$ for all $(t,x)\in [0,T]\times\R^d$, and let $\Pi^g_c$ be the subset of jointly continuous functions.
  \item Given a filtration $\bbG$, we let $\Pred(\bbG)$ denote the $\sigma$-algebra of $\bbG$-predictably measurable subsets of $[0,T]\times \Omega$.
  \item We let $\mcT$ be the set of all $[0,T]$-valued $\bbF$-stopping times and for each $\eta\in\mcT$, we let $\mcT_\eta$ be the subset of stopping times $\tau$ such that $\tau\geq \eta$, $\Prob$-a.s. Moreover, we let $\mcT^t$ (resp. $\mcT^t_\eta$) be the corresponding subsets of $\bbF^t$-stopping times, with $\tau\geq t$ (resp. $\tau\geq \eta$), $\Prob$-a.s.
  \item For $p\geq 1$, $t\in [0,T]$ and $\tau\in\mcT^t$, we let $\mcS^{p}_{[t,\tau]}$ be the set of all $\R$-valued, $\bbF^t$-progressively measurable, \cadlag processes $(Z_s: s\in [t,\tau])$ for which $\|Z\|_{\mcS^p_{[t,\tau]}}:=\E\big[\sup_{s\in[t,\tau]} |Z_s|^p\big]<\infty$. Moreover, we let $\mcA^p_{[t,\tau]}$ be the subset of $\bbF^t$-predictably measurable and non-decreasing processes with $Z_t=0$. Whenever $\tau=T$ we use the notations $\mcS^p_{t}$ and $\mcA^p_{t}$, respectively.
  \item We let $\mcH^{p}_{[t,\tau]}(W)$ denote the set of all $\R^d$-valued $\mcP(\bbF^t)$-measurable processes $(Z_s: s\in[t,\tau])$ such that $\|Z\|_{\mcH^p_{[t,\tau]}(W)}:=\E\big[\big(\int_t^\tau |Z_s|^2 ds\big)^{p/2}\big]^{1/p}<\infty$. Furthermore, we set $\mcH^{p}_{t}(W):=\mcH^{p}_{[t,T]}(W)$.
  \item We let $\mcH^{p}_{[t,\tau]}(\mu)$ denote the set of all $\R$-valued, $\mcP(\bbF^t) \otimes\mcB(E)$-measurable maps $(Z_s(e): s\in[t,\tau],e\in E)$ such that $\|Z\|_{\mcH^p_{[t,\tau]}(\mu)}:=\E\big[\big(\int_t^\tau \!\!\int_E |Z_s(e)|^2 \lambda(de)ds\big)^{p/2}\big]^{1/p}<\infty$ and set $\mcH^{p}_{t}(\mu)=\mcH^{p}_{[t,T]}(\mu)$.
  \item We let $\bigS^{p}$ be the set of all maps $Z:\cup_{t\in[0,T]}[t,T]\times \Omega\times \{t\}\times \R^d\to \R:(s,\omega,t,x)\mapsto Z^{t,x}_s(\omega)$ such that $Z^{t,x}\in \mcS^p_t$ and there is a $v\in\Pi^g_c$ such that $v(t,x)=Z^{t,x}_t$, $\Prob$-a.s., for all $(t,x)\in [0,T]\times\R^n$. Moreover, given $Z\in\bigS^p$, we let $(t,x)\mapsto Z(t,x)$ denote this deterministic function, so that $ Z(t,x)=v(t,x)$ for all $(t,x)\in[0,T]\times \R^d$.
  \item For $0\leq t\leq s\leq T$, we let $\bar \mcU_s$ be the set of all $u=(\tau_j,\beta_j)_{1\leq j\leq N}$, where $(\tau_j)_{j=1}^\infty$ is a non-decreasing sequence of $\bbF$-stopping times in $\mcT_s$, $\beta_j$ is a $\mcF_{\tau_j}$-measurable r.v.~taking values in the compact set $U$ and $N:=\max\{j:\tau_j<T\}$, such that $\Xi^{t,x;u}_T\in L^2(\Omega,\mcF_T,\Prob)$, where $\Xi^{t,x;u}_T$ is the total cost of impulses for the non-randomizing player (see \eqref{ekv:Xi-def}).
  \item For each $k\in\bbN$, we let $\bar\mcU^{k}_s$ be the subset of $\bar\mcU_s$ consisting of all $u=(\tau_j,\beta_j)_{1\leq j\leq N}$ for which $N\leq k$, $\Prob$-a.s.
  \item For $u\in\bar\mcU$, we let $[u]_{j}:=(\tau_i,\beta_i)_{1\leq i\leq N\wedge j}$. Moreover, we introduce $N(s):=\max\{j\geq 0:\tau_j\leq s\}$ and let $u_s:=[u]_{N(s)}$.
  \item We let $\bar\mcV$ denote the set of all $\Pred(\bbF)\otimes\mcB(E)$-measurable, bounded maps $\nu:[0,T]\times\Omega\times E\to [0,\infty)$ and for each $n\in\bbN$, we denote by $\bar\mcV^n$ the subset of maps $\nu:[0,T]\times\Omega\times E\to [0,n]$.
\end{itemize}

Next, we define the set of \emph{feedback impulse controls}.
\begin{defn}
For $(t,x)\in [0,T]\times\R^d$ and $\tilde u\in\bar\mcU$, we let $\mcG^{t,x,\tilde u}$ be the augmented natural filtration generated by $X^{t,x,\tilde u}$. For $s\in [t,T]$ and $k\in\bbN$, we then denote by $\mcU^{t,x}_s$ (resp.~$\mcU^{t,x,k}_s$) the subset of all $u\in\bar\mcU_s$ (resp.~$u\in\bar\mcU^k_s$) such that $\tau_j$ is a $(\mcG^{t,x,[u]_{j-1}}_r)_{t\leq r\leq T}$-stopping time and $\beta_j$ is $\mcG^{t,x,[u]_{j-1}}_{\tau_j}$-measurable.
\end{defn}

The set of \emph{feedback randomization densities} is defined similarly:
\begin{defn}
We denote by $\mcV^{t,x}$ the set of all maps $\bar \mcU\to\bar\mcV$ ($u\mapsto\nu(u)$) such that $\nu(u)$ is a $\Pred(\bbG^{t,x,u})\otimes\mcB(E)$-measurable, bounded map $\nu(u):[0,T]\times\Omega\times E\to [0,\infty)$ and for each $n\in\bbN$, we denote by $\mcV^{t,x,n}$ the subset of $\nu\in\mcV^{t,x}$, containing all maps $\nu:[0,T]\times\Omega\times E\to [0,n]$. Whenever it is apparent from the context which impulse control that $u$ refers to, we simplify the notation by writing $\nu$ instead of $\nu(u)$.
\end{defn}

%\begin{defn}
%For $t\in[0,T]$, the set of \emph{non-anticipative impulse strategies} $\mcU^{t,S}_s$ is defined as the set of all maps $u^S:\mcV\to \mcU^t_s$ such that for any $\nu,\tilde \nu\in\mcV$, we have $[u^S(\tilde\nu)]_s=[u^S(\nu)]_s$, whenever $\tilde\nu(r,\omega,e)=\nu(r,\omega,e)$, $dr\times de\times d\Prob$-a.e.~on $[t,s]\times\Omega\times E$. Moreover, we let $\mcU^{k,S}_t$ be the subset of all maps $u^S:\mcV\to \mcU^{k}_t$.
%\end{defn}
%
%
%\begin{defn}
%The set of \emph{non-anticipative density strategies} $\mcV^S_t$ is defined as the set of all maps $\nu^S:\mcU\to\mcV_t$ such that for any $u,\tilde u\in\mcU$, $\ett_{[[\tilde u]_s=[u]_s]}\nu^S(\tilde u)(s,e)=\ett_{[[\tilde u]_s=[u]_s]}\nu^S(\tilde u)(s,e)$, $ds\times \lambda(de)\times d\Prob$-a.e.~and $\nu^S(u)$ is bounded, uniformly in $u$.% Moreover, we let $\mcV^{n,S}_t$ be the subset of all maps $u^S: \mcU_t\to\mcV^n$.
%\end{defn}

We also mention that, unless otherwise specified, all inequalities involving random variables are assumed to hold $\Prob$-a.s.

\subsection{Assumptions}
\begin{ass}\label{ass:lambda}
 $\lambda$ is a finite positive measure on $(U,\mcB(U))$ with full topological support.
\end{ass}

We assume that the coefficients of the forward SDE satisfy the following constraints:
\begin{ass}\label{ass:onSDE}
%For any $t \in [0,T]$, $b\in U$, $e\in E$ and $x,x'\in\R^d$ we have:
\begin{enumerate}[i)]
  \item\label{ass:onSDE-Gamma} The functions $\xi:[0,T]\times\R^d\times U\to\R^d$ and $\gamma:[0,T]\times\R^d\times E\to\R^d$ are jointly continuous and satisfy the growth condition
  \begin{align}\label{ekv:imp-bound}
   |x+\xi(t,x,b)| \vee |x+\gamma(t,x,e)|\leq K_{\gamma,\xi}\vee |x|
  \end{align}
  for some constant $K_{\gamma,\xi}>0$. Moreover, $\gamma$ is Lipschitz continuous in $x$ uniformly in $(t,e)$, \ie
  \begin{align*}
    |\gamma(t,x',e)-\gamma(t,x,e)|\leq k_\gamma|x'-x|
  \end{align*}
  for all $(t,x,x',e)\in [0,T]\times \R^d\times\R^d\times E$.
  \item\label{ass:onSDE-a-sigma} The coefficients $a:[0,T]\times\R^d\to\R^{d}$ and $\sigma:[0,T]\times\R^d\to\R^{d\times d}$ satisfy the growth condition
  \begin{align*}
    |a(t,x)|+|\sigma(t,x)|&\leq C(1+|x|)
  \end{align*}
  and the Lipschitz continuity
  \begin{align*}
    |a(t,x)-a(t,x')|+|\sigma(t,x)-\sigma(t,x')|&\leq k_{a,\sigma}|x'-x|
  \end{align*}
  for all $(t,x,x')\in [0,T]\times \R^d\times\R^d$.
\end{enumerate}
\end{ass}

Moreover, we make the following assumption on the coefficients of the BSDE:

\begin{ass}\label{ass:oncoeff}
For some $\rho\geq 1$ we have:
\begin{enumerate}[i)]
  \item\label{ass:oncoeff-f} The running cost $f:[0,T]\times \R^{d+1+d}\to\R$ is jointly continuous in $(t,x)$, uniformly in $(y,z)$, satisfies the growth condition
  \begin{align*}
    |f(t,x,0,0)|\leq C(1+|x|^\rho)
  \end{align*}
  and the Lipschitz condition
  \begin{align*}
    |f(t,x,y,z)-f(t,x,y',z')|&\leq k_f(|y-y|+|z-z'|)
  \end{align*}
  for all $(t,x,y,y',z,z')\in [0,T]\times \R^d\times\R\times\R\times\R^d\times\R^d$.
  \item The terminal value $\psi:\R^d\to\R$ is continuous, of polynomial growth, \ie
  \begin{align*}
    |\psi(x)|\leq C(1+|x|^\rho)
  \end{align*}
  and satisfies
  \begin{align*}
    \sup_{b\in U}\big\{\psi(x+\xi(T,x,b))-\ell(T,x,b)\big\}\leq \psi(x) \leq \inf_{e\in E}\big\{\psi(x+\gamma(T,x,e))+\chi(T,x,e)\big\}
  \end{align*}
  for all $x\in \R^d$.
  \item\label{ass:oncoeff-ell} The intervention costs $\ell:[0,T]\times\R^d\times U\to [\delta,\infty)\subset (0,\infty)$ and $\chi:[0,T]\times\R^d\times E\to[0,\infty)$ are jointly continuous and of polynomial growth.
\end{enumerate}
\end{ass}

Moreover, we let $\supOP v(T,x)=\infOP v(T,x)=v(T,x)$. We introduce the following assumption, which is used solely to prove the comparison result for viscosity solutions to \eqref{ekv:var-ineq}:

\begin{ass}\label{ass:nofreeloop}
For each $R>0$, there are functions $h_1,h_2:[0,T]\to (0,\infty)$ such that for each $(t,x_0)\in [0,T]\times \bar B_R(\R^d)$ it holds that, whenever $(b_i)_{i\in\bbN}$, $(e_i)_{i\in\bbN}$ and $(\iota_i)_{i\in \bbN}$ are sequences in $U$, $E$ and $\{1,2\}$, respectively, for which there exists a $\kappa\in\bbN$ such that $|x_{\kappa}-x_0|\leq h_1(t)$, where
\begin{align*}
  x_{j}=x_{j-1}+\ett_{[\iota_j=1]}\xi(t,x_{j-1},b_j)+\ett_{[\iota_j=2]}\gamma(t,x_{j-1},e_j),\quad \text{for }j=1,2,\ldots,\kappa
\end{align*}
then
\begin{align*}
  |\sum_{j=1}^{\kappa}(\ett_{[\iota_j=1]}\ell(t,x_{j-1},b_j)-\ett_{[\iota_j=2]}\chi(t,x_{j-1},e_j))|\geq h_2(t).
\end{align*}
\end{ass}

The assumption above is a variant of the no-free-loop condition commonly used in optimal switching problems (see \eg assumption H3.(ii) in \cite{Morlais13} or assumption H4 in \cite{BollanSWG1}). It is intended to prevent the players from becoming trapped in an infinite cycle of simultaneous impulses.

\subsection{Viscosity solutions}

We define the upper, $v^*$, and lower, $v_*$, semi-continuous envelope of a function $v:[0,T]\times\R^d\to\R$ as
\begin{align*}
v^*(t,x):=\limsup_{(t',x')\to(t,x),\,t'<T}v(t',x')\quad {\rm and}\quad v_*(t,x):=\liminf_{(t',x')\to(t,x),\,t'<T}v(t',x').
\end{align*}
Next, we introduce the limiting parabolic superjet $\bar J^+v$ and subjet $\bar J^-v$.
\begin{defn}\label{def:jets}
Subjets and superjets
\begin{enumerate}[i)]
\item For a l.s.c. (resp. u.s.c.) function $v : [0, T]\times \R^d \to \R$, the parabolic subjet, denote by $J^-v(t, x)$, (resp.~the parabolic superjet, $J^+v(t, x)$) of $v$ at $(t, x) \in [0, T]\times \R^d$, is defined as the set of triples $(p, q,M) \in \R \times\R^d \times \mathbb S^d$ satisfying
    \begin{align*}
      v(t', x') \geq (\text{resp.~}\leq)\: v(t, x) + p(t'- t)+ < q, x' - x > + \tfrac{1}{2} < x' - x,M(x' - x) > +o(|t' - t| + |x' - x|^2)
    \end{align*}
    for all $(t',x')\in[0,T]\times\R^d$, where $\mathbb S^d$ is the set of symmetric real matrices of dimension $d\times d$.
\item For a l.s.c. (resp. u.s.c.) function $v : [0, T]\times \R^d \to \R$ we denote by $\bar J^-v(t, x)$ the parabolic limiting
subjet (resp. $\bar J^+v(t, x)$ the parabolic limiting superjet) of $v$ at $(t, x) \in [0, T]\times \R^d$, defined as the set of triples $(p, q,M) \in \R \times\R^d \times \mathbb S^d$ such that:
\begin{align*}
  (p, q,M) = \lim_{n\to\infty} (p_n, q_n,M_n),\quad (t, x) = \lim_{n\to\infty}(t_n, x_n)
\end{align*}
for some sequence $(t_n,x_n,p_n,q_n,M_n)_{n\geq 1}$ with $(p_n, q_n,M_n) \in J^-v(t_n, x_n)$ (resp. $(p_n, q_n,M_n) \in J^+v(t_n, x_n)$) for all $n\geq 1$ and $v(t, x) = \lim_{n\to\infty}v(t_n, x_n)$.
\end{enumerate}
\end{defn}

We now give the definition of a viscosity solution for the QVI in \eqref{ekv:var-ineq}. (see also pp. 9-10 of \cite{UsersGuide}).
\begin{defn}\label{def:visc-sol-jets}
Let $v$ be a locally bounded function from $[0,T]\times \R^d$ to $\R$. Then,
\begin{enumerate}[a)]
  \item It is referred to as a viscosity supersolution (resp. subsolution) to \eqref{ekv:var-ineq} if it is l.s.c.~(resp. u.s.c.) and satisfies:
  \begin{enumerate}[i)]
    \item $v(T,x)\geq \psi(x)$ (resp. $v(T,x)\leq \psi(x)$)
    \item For any $(t,x)\in [0,T)\times\R^d$ and $(p,q,X)\in \bar J^- v(t,x)$ (resp. $\bar J^+ v(t,x)$) we have
    \begin{align*}
      \min\big\{&v(t,x)-\supOP v(t,x)(t,x),\max\{v(t,x)-\infOP v(t,x),-p-q^\top a(t,x)
      \\
      &-\frac{1}{2}\trace(\sigma\sigma^\top(t,x)X)-f(t,x,v(t,x),\sigma^\top(t,x)q)\}\big\}\geq 0\quad (\text{resp. }\leq 0)
    \end{align*}
  \end{enumerate}
  \item It is called a viscosity solution to \eqref{ekv:var-ineq} if $v_*$ is a supersolution and $v^*$ is a subsolution.
\end{enumerate}
\end{defn}

We will sometimes use the following alternative definition of viscosity supersolutions (resp. subsolutions):
\begin{defn}\label{def:visc-sol-dom}
  A l.s.c.~(resp. u.s.c.) function $v$ is a viscosity supersolution (resp.~subsolution) to \eqref{ekv:var-ineq} if $v(T,x)\leq \psi(x)$ (resp. $\geq \psi(x)$) and whenever $\varphi\in C^{1,2}([0,T]\times\R^d\to\R)$ is such that $\varphi(t,x)=v(t,x)$ and $\varphi-v$ has a local maximum (resp. minimum) at $(t,x)$, then
  \begin{align*}
    \min\big\{&v(t,x)-\supOP v(t,x)(t,x),\max\{v(t,x)-\infOP v(t,x),-\varphi_t(t,x)-\mcL\varphi(t,x)
    \\
    &-f(t,x,v(t,x),\sigma^\top(t,x)\nabla_x\varphi(t,x))\}\big\}\geq 0\quad(\text{resp. }\leq 0).
  \end{align*}
\end{defn}

\begin{rem}\label{rem:mcM-monotone}
Let $u,v:[0,T]\times \R^d\to\R$ be locally bounded functions. We remark that $\supOP$ and $\infOP$ are both monotone (if $u\leq v$ pointwise, then $\supOP u\leq \supOP v$ and $\infOP u\leq \infOP v$). Moreover, $\supOP(u_*)$ and $\infOP(u_*)$ (resp. $\supOP(u^*)$ and $\infOP(u^*)$) are l.s.c.~(resp. u.s.c.).% and $(\supOP u)_*\leq

In particular, $\supOP v$ and $\infOP v$ are both jointly continuous whenever $v$ is.\\
\end{rem}

\subsection{Some preliminary result}

Some of the results gathered in this section uses the auxiliary lower barrier $h\in \Pi^g_c$ which is assumed to satisfy $h(T,\cdot)\leq\psi$. Hence, $\Psi(t,x):=\ett_{[t<T]}h(t,x)+\ett_{[t=T]}\psi(x)$ is jointly continuous on $[0,T)\times\R^d$ and $\lim_{(t,x')\to(T,x)}\Psi(t,x')\leq \Psi(T,x)$ for all $x\in\R^d$.

We give the following proposition, strategically formulated to streamline subsequent implementation processes. A proof (based on comparable findings documented in \cite{Dumitrescu15} and \cite{HamMor16}) can be found in \cite{qvi-stop}.

\begin{prop}\label{prop:rbsde-jmp}
For each $(t,x)\in [0,T]\times\R^d$ and $n\in\bbN$, there is a unique quadruple\\ $(Y^{t,x,n},Z^{t,x,n},V^{t,x,n},K^{+,t,x,n})\in\mcS^{2}_t \times \mcH^{2}_t(W) \times \mcH^{2}_t(\mu) \times \mcA^{2}_{t}$ such that
\begin{align}\label{ekv:rbsde-jmp-base}
  \begin{cases}
     Y^{t,x,n}_s=\psi(X^{t,x}_T)+\int_s^T f^n(r,X^{t,x}_r, Y^{t,x,n}_r, Z^{t,x,n}_r, V^{t,x,n}_{r})dr -\int_s^T  Z^{t,x,n}_r dW_r
    \\
    \quad-\int_{s}^T\!\!\!\int_E  V^{t,x,n}_{r}(e)\mu(dr,de)+(K^{+,t,x,n}_T- K^{+,t,x,n}_s),\quad \forall s\in [t,T],
    \\
    Y^{t,x,n}_s\geq h(s,X^{t,x}_s),\, \forall s\in [t,T] \quad\text{and}\quad\int_t^T\!\! \big(Y^{t,x,n}_s-h(s,X^{t,x}_s)\big)dK^{+,t,x,n}_s,
  \end{cases}
\end{align}
where
\begin{align*}
f^n(t,x,y,z,v):= f(t,x,y,z)-n\int_E(v(e)+\chi(t,x,e))^-\lambda(de).
\end{align*}
Moreover, $Y^{t,x,n}_s=\esssup_{\tau\in\mcT_s}Y^{t,x,n;\tau}_s$, where for each $\tau\in\mcT_t$, the triple
$(Y^{t,x,n;\tau},Z^{t,x,n;\tau},V^{t,x,n;\tau})\in\mcS^{2}_{[0,\tau]} \times \mcH^{2}_{[0,\tau]}(W) \times \mcH^{2}_{[0,\tau]}(\mu)$ satisfies
\begin{align}\nonumber
  Y^{t,x,n;\tau}_s&=\Psi(\tau,X^{t,x}_\tau)+\int_s^\tau f^n(r,X^{t,x}_r, Y^{t,x,n;\tau}_r, Z^{t,x,n;\tau}_r, V^{t,x,n;\tau}_{r})dr
  \\
  &\quad -\int_s^\tau  Z^{t,x,n;\tau}_r dW_r-\int_{s}^\tau\!\!\!\int_E  V^{t,x,n;\tau}_{r}(e)\mu(dr,de),\quad \forall s\in [0,\tau],\label{ekv:bsde-jmp}
\end{align}
and for each $\eta\in \mcT_t$, the stopping time
\begin{align*}
  \tau^*:=\inf\{s\geq \eta:Y^{t,x,n}_s=h(s,X^{t,x}_s)\}\wedge T
\end{align*}
is optimal in the sense that $Y^{t,x,n}_\eta=Y^{t,x,n;\tau^*}_\eta$. Finally, there is a function $v_n\in\Pi^g_c$ such that $v_n(s,X^{t,x}_s)=Y^{t,x,n}_s$ for all $s\in[t,T]$ and $v_n$ is the unique viscosity solution (see Definition~\ref{rem:visc-non-loc} below for a definition) in $\Pi^g_c$ to
\begin{align}\label{ekv:obst-prob-n}
\begin{cases}
  \min\{v_n(t,x)-h(t,x),-\frac{\partial}{\partial t}v_n(t,x)-\mcL v_n(t,x)+\mcK^n v_n(t,x)\\
  \quad- f(t,x,v_n(t,x),\sigma^\top(t,x)\nabla_x v_n(t,x))\}=0,\quad\forall (t,x)\in[0,T)\times \R^d \\
  v_n(T,x)=\psi(x),
\end{cases}
\end{align}
where $\mcK^n \phi(t,x):=n\int_E(\phi(t,x+\gamma(t,x,e)) + \chi(t,x,e)-\phi(t,x))^-\lambda(de)$.
\end{prop}

Note that the variational inequality in \eqref{ekv:obst-prob-n} is non-local in the sense that, by invoking the operator $\mcK^n$, the driver depends on values of $v$ at points other that $(t,x)$. When defining the corresponding viscosity solutions, we adhere to the principle in \cite{HamMor16} (as opposed to the path taken in \cite{Barles97}) and define:

\begin{defn}\label{rem:visc-non-loc}
A locally bounded map $v:[0,T]\times\R^d\to\R$ is a viscosity supersolution (resp. subsolution) to \eqref{ekv:obst-prob-n} if it is l.s.c.~(resp. u.s.c.) satisfies $v(T,x)\geq \psi(x)$ (resp. $v(T,x)\leq \psi(x)$) and if for any $(t,x)\in [0,T)\times\R^d$ and $(p,q,X)\in \bar J^- v(t,x)$ (resp. $\bar J^+ v(t,x)$) we have
\begin{align*}
  \min\big\{&v(t,x)-h(t,x),-p-q^\top a(t,x)
  \\
  &-\frac{1}{2}\trace(\sigma\sigma^\top(t,x)X)+\mcK^n v(t,x)+f(t,x,v(t,x),\sigma^\top(t,x)q)\}\big\}\geq 0\quad (\text{resp. }\leq 0).
\end{align*}
It is a viscosity solution if $v_*$ is a supersolution and $v^*$ is a subsolution.
\end{defn}

\begin{rem}\label{rem:Vtx-cont}
The proof of Proposition~\ref{prop:rbsde-jmp} delineated in \cite{qvi-stop} borrows an important fact from \cite{HamMor16} (see Proposition 3.1 therein), namely that the jump component can be written $V^{t,x,n}_s(e)=v_n(s,X^{t,x}_s+\gamma(s,X^{t,x}_s,e))-v_n(s,X^{t,x}_s)$. Of particular importance for our subsequent analysis is the fact that the map $(t,x,e)\mapsto V^{t,x,n}_t(e)$ can be choosen to be deterministic and continuous.
\end{rem}

\bigskip

For each $(t,x)\in [0,T]\times\R^d$, the sequence $(Y^{t,x,n})_{n\in\bbN}$ is a non-increasing sequence of \cadlag processes that are bounded from below by the process $\Psi(\cdot,X^{t,x}_\cdot)$. There is thus a progressively measurable process $Y^{t,x}$ such that $Y^{t,x,n}\searrow Y^{t,x}$ pointwisely, as $n\to\infty$. The following result, which is a main result of \cite{qvi-stop} (see Theorem 3.2 and Proposition 3.3 therein and also Theorem 3.1 in \cite{imp-stop-game}), provides a characterization of $Y^{t,x}$ in the form of a non-linear Snell envelope:

\begin{thm}\label{thm:qvi-stop}
For each $(t,x)\in [0,T]\times\R^d$, the process $Y^{t,x}$ belongs to $\mcS^2_t$ and satisfies $Y^{t,x}_\eta=\esssup_{\tau\in\mcT_\eta}Y^{t,x;\tau}_\eta$ for each $\eta\in\mcT_t$, where $(Y^{t,x;\tau},Z^{t,x;\tau},V^{t,x;\tau},K^{-,t,x;\tau})\in\mcS^{2}_{[0,\tau]} \times \mcH^{2}_{[0,\tau]}(W) \times \mcH^{2}_{[0,\tau]}(\mu) \times \mcA^{2}_{[0,\tau]}$ is the unique maximal solution to the BSDE
\begin{align}\label{ekv:stopped-bsde}
  \begin{cases}
     Y^{t,x;\tau}_s=\Psi(\tau,X^{t,x}_\tau)+\int_s^\tau f(r,X^{t,x}_r, Y^{t,x;\tau}_r, Z^{t,x;\tau}_r)dr -\int_s^\tau  Z^{t,x;\tau}_r dW_r-\int_{s}^\tau\!\!\!\int_E  V^{t,x;\tau}_{r}(e)\mu(dr,de)
    \\
    \quad-(K^{-,t,x;\tau}_\tau- K^{-,t,x;\tau}_s),\quad \forall s\in [0,\tau],
    \\
    V^{t,x;\tau}_s(e)\geq - \chi(s,X^{t,x}_{s-},e),\quad d\Prob\otimes ds\otimes \lambda(de)-a.e.,
  \end{cases}
\end{align}
Moreover, for each $\eta\in \mcT_t$, the stopping time
\begin{align*}
  \tau^*:=\inf\{s\geq \eta:Y^{t,x}_s=\Psi(s,X^{t,x}_s)\}
\end{align*}
is optimal in the sense that $Y^{t,x}_\eta=Y^{t,x;\tau^*}_\eta$. Finally, $v(t,x)=Y^{t,x}_t$ belongs to $\Pi^g_c$ and is the unique viscosity solution to
\begin{align}\label{ekv:dbl-obst-pde}
\begin{cases}
  \min\{v(t,x)-h(t,x),\max\{v(t,x)-\infOP v(t,x),-v_t(t,x)-\mcL v(t,x)\\
  \quad-f(t,x,v(t,x),\sigma^\top(t,x)\nabla_x v(t,x))\}\}=0,\quad\forall (t,x)\in[0,T)\times \R^d \\
  v(T,x)=\psi(x).
\end{cases}
\end{align}
Here, viscosity solutions to \eqref{ekv:dbl-obst-pde} are defined analogously to viscosity solutions of \eqref{ekv:var-ineq}.
\end{thm}

Whereas the unique solution to \eqref{ekv:dbl-obst-pde}, for a suitably chosen non-decreasing sequence of maps $h$, will provide a means of approximating a viscosity supersolution to \eqref{ekv:var-ineq} from below, the following proposition will allow us to approximate a subsolution from above:

\begin{prop}\label{prop:qvi-rbsde}
There is a unique viscosity solution in $\Pi^g_c$ to the following quasi-variational inequality
\begin{align}\label{ekv:var-ineq-n}
\begin{cases}
  \min\{v(t,x)-\supOP v(t,x),-v_t(t,x)-\mcL v(t,x)+\mcK^n v(t,x)\\
  \quad-f(t,x,v(t,x),\sigma^\top(t,x)\nabla_x v(t,x))\}=0,\quad\forall (t,x)\in[0,T)\times \R^d \\
  v(T,x)=\psi(x).
\end{cases}
\end{align}
\end{prop}

\noindent\emph{Proof.} The result is a special case of Theorem 4.5 in \cite{qvi-rbsde}.\qed\\

%%%%%%%%%%%%%%%%%%%%%%%%%%%%%%%%%%%%%%%%%%%%%%%%%%%%%%%%%%%
%%%%%%%%%%%%%%%%%%%%%%%%%%%%%%%%%%%%%%%%%%%%%%%%%%%%%%%%%%%
%%%%%%%%%%%%%%%%%%%%%%%%%%%%%%%%%%%%%%%%%%%%%%%%%%%%%%%%%%%

\section{A game of impulse control with randomized controls\label{sec:game}}

For $\nu\in\bar\mcV$ and $u\in\bar\mcU_t$, we let $(P^{t,x;u,\nu},Q^{t,x;u,\nu},R^{t,x;u,\nu})$ be the unique triple with $P^{t,x;u,\nu}-\Xi^{t,x;u}\in\mcS^2$ and $(Q^{t,x;u,\nu},R^{t,x;u,\nu})\in\mcH^2(W)\times\mcH^2(\mu)$, that satisfies
\begin{align}\nonumber
P^{t,x;u,\nu}_s&=\psi(X^{t,x;u}_T)+\int_s^T f^\nu(r,X^{t,x;u}_r,P^{t,x;u,\nu}_r,Q^{t,x;u,\nu}_r,R^{t,x;u,\nu}_{r})dr-\int_s^T Q^{t,x;u,\nu}_r dW_r
\\
&\quad-\int_{s}^T\!\!\!\int_E R^{t,x;u,\nu}_{r}(e)\mu(dr,de) - \Xi^{t,x;u}_{T}+\Xi^{t,x;u}_s, \label{ekv:non-ref-bsde-simp}
\end{align}
with driver
\begin{align*}
f^\nu(t,x,y,z,v):=f(t,x,y,z)+\int_E(v(e)+\chi(t,x,e))\nu_t(e)\lambda(de).
\end{align*}
Note that the BSDEs in \eqref{ekv:non-ref-bsde-simp} and \eqref{ekv:non-ref-bsde-simp-alt} are equivalent in the sense that $(P^{t,x;u,\nu},Q^{t,x;u,\nu},R^{t,x;u,\nu})$ is also a solution to \eqref{ekv:non-ref-bsde-simp-alt}.

The main results of this work is summarized in the following theorem:

\begin{thm}\label{thm:game-value}
The zero-sum game of control versus randomized control with reward/cost $P^{t,x;\cdot,\cdot}_s$ has a value and for each $(t,x)\in[0,T]\times\R^d$ and $s\in[t,T]$, we have
\begin{align}\label{ekv:game-value}
  \essinf_{\nu\in \mcV^{t,x}}\esssup_{u\in\mcU^{t,x}_{s}}P^{t,x;u,\nu}_{s}=\esssup_{u\in\mcU^{t,x}_{s}}\essinf_{\nu\in\mcV^{t,x}}P^{t,x;u,\nu}_{s} = v(s,X^{t,x}_{s}),
\end{align}
where $v\in\Pi^g_c$ is the unique viscosity solution to \eqref{ekv:var-ineq} within $\Pi^g$. Furthermore, we have the representation $v({s},X^{t,x}_{s})=Y^{t,x}_{s}$, where $Y$ is the unique element in $\bigS^2$ that satisfies \eqref{ekv:Snell-env}-\eqref{ekv:fwd-sde}.
\end{thm}

\subsection{Some preliminary estimates}
Before we proceed, we give some preliminary estimates that will be useful later on.

For $\nu\in\mcV$, we let $\E^\nu$ be expectation with respect to the probability measure $\Prob^\nu$ on $(\Omega,\mcF)$ defined by $d\Prob^\nu:=\kappa^\nu_T d\Prob$. Here,
\begin{align*}
\kappa^{\nu}_s&:=\mcE_s\Big(\int_{0}^\cdot\!\!\int_E(\nu_r(e)-1)(\mu(dr,de)-\lambda(de)dr)\Big)
\\
&:=\exp\Big(\int_{0}^s\!\!\!\int_E(1-\nu_r(e))\lambda(de)dr\Big)\prod_{\sigma_j\leq s}\nu_{\sigma_j}(\zeta_j)
\end{align*}
and the sequence $(\sigma_j,\zeta_j)_{j\geq 1}$ is the one that appears in the Dirac sum formulation of the random measure, $\mu=\sum_{j\geq 1}\delta_{(\sigma_j,\zeta_j)}$.

\begin{prop}\label{prop:SDEmoment}
%Under Assumption~\ref{ass:on-coeff},
For each $p\geq 1$, there is a $C>0$ such that
\begin{align}\label{ekv:SDEmoment}
\E^\nu\Big[\sup_{s\in[\eta,T]}|X^{t,x;u}_s|^{p}\Big|\mcF_\eta\Big]\leq C(1+|X^{t,x;u}_\eta|^{p}),
\end{align}
$\Prob$-a.s.~for all $(t,x)\in [0,T]\times\R^d$ and $(\eta,u,\nu)\in \mcT_t\times\mcU_t\times\mcV_t$.
\end{prop}

\noindent\emph{Proof.} Let $\tilde \mcU$, be the set of all $\tilde u:=(\theta_j,\vartheta_j)_{j=1}^{M}$ where $(\theta_j)_{j=1}^\infty$ is an increasing sequence of $\bbF$-stopping times with $\theta_1\geq \eta$, the random variable $\vartheta_j$ is $\mcF_{\theta_j}$-measurable and $E$-valued and $M:=\inf\{j\geq 0:\theta_j\leq T\}$. If $X^{t,x;u,v}$ solves the SDE
\begin{align*}
X^{t,x;u,\tilde u}_s&=x+\int_t^s a(r,X^{t,x;u,\tilde u}_r)dr+\int_t^s\sigma(r,X^{t,x;u}_r)dW_r + \int_t^{s\wedge\eta}\!\!\!\int_E\gamma(r,X^{t,x;u,\tilde u}_{r-},e)\mu(dr,de)
\\
&\quad + \sum_{j=1}^{M}\ett_{[\theta_j\leq s]}\gamma(\theta_j,X^{t,x;u,\tilde u}_{\theta_j-},\vartheta_j)+\sum_{j=1}^{N}\ett_{[\tau_j\leq s]}\xi(\tau_j,X_{\tau_j}^{t,x;[u]_{j-1}},\beta_j),\quad\forall s\in[t,T].
\end{align*}
Then a standard argument in control randomization (see \eg Section 4.1 in \cite{Bandini18}) implies that
\begin{align*}
\esssup_{u\in\mcU_t}\esssup_{\nu\in\mcV_t}\E^\nu\Big[\sup_{s\in[\eta,T]}|X^{t,x;u}_s|^{p}\Big|\mcF_\eta\Big]\leq \esssup_{u\in\mcU_t}\esssup_{\tilde u\in\tilde\mcU_t}\E\Big[\sup_{s\in[\eta,T]}|X^{t,x;u,\tilde u}_s|^{p}\Big|\mcF_\eta\Big].
\end{align*}
On the other hand, by arguing as in the proof of Proposition 5.2 in~\cite{Stochastics23} we find that there is a $C>0$ such that
\begin{align*}
\esssup_{u\in\mcU_t}\esssup_{\tilde u\in\tilde\mcU_t}\E\Big[\sup_{s\in[\eta,T]}|X^{t,x;u,\tilde u}_s|^{p}\Big|\mcF_\eta\Big]\leq C(1+|X^{t,x;u}_\eta|^{p})
\end{align*}
$\Prob$-a.s.~for all $(t,x)\in [0,T]\times\R^d$ and $(\eta,u,\tilde u)\in \mcT_t\times\mcU_t\times \tilde \mcU_t$.\qed\\

For $(t,x)\in[0,T]\times\R^d$, we introduce the sequence $(Y^{t,x,0,n})_{n\in\bbN}$, by for each $n\in\bbN$, letting the process $Y^{t,x,0,n}\in\mcS^2_t$ be the first component in the unique solution to the standard BSDE with jumps (recall the definition of $X^{t,x}$ in \eqref{ekv:fwd-sde})
\begin{align}\nonumber
    Y^{t,x,0,n}_s&=\psi(X^{t,x}_T)+\int_s^T f(r,X^{t,x}_r,Y^{t,x,0,n}_r,Z^{t,x,0,n}_r)dr-n\int_t^T\!\!\int_E(V^{t,x,0,n}_r(e)+\chi(r,X^{t,x}_{r},e))^-\lambda(de)dr
    \\
    &\quad -\int_s^T Z^{t,x,0,n}_r dW_r-\int_{s}^T\!\!\!\int_E V^{t,x,0,n}_{r}(e)\mu(dr,de),\quad \forall s\in [t,T],\label{ekv:rbsde-pen-0}
\end{align}
Using the comparison principle for BSDEs with jumps, we easily deduce the following moment estimates:
\begin{prop}\label{prop:BSDEmoment}
%Under Assumption~\ref{ass:on-coeff},
There is a $C>0$ such that for each $n\in\bbN$, we have
\begin{align}\label{ekv:Ybnd}
  \|Y^{t,x,0,n}\|_{\mcS^2_t}\leq C(1+|x|^\rho)
\end{align}
for all $(t,x)\in [0,T]\times\R^d$.
\end{prop}

\noindent\emph{Proof.} First, note that $Y^{t,x,0,n}\leq Y^{t,x,0,0}$, where the latter satisfies $\|Y^{t,x,0,0}\|_{\mcS^2_t}\leq C(1+|x|^\rho)$. On the other hand, $Y^{t,x,0,n}\geq \tilde Y^{t,x,n}$, where $(\tilde Y^{t,x,n},\tilde Z^{t,x,n},\tilde V^{t,x,n})\in\mcS^2_t\times\mcH^2_t(W)\times\mcH^2_t(\mu)$ satisfies
\begin{align*}
    \tilde Y^{t,x,n}_s&=\psi(X^{t,x}_T)+\int_s^T f(r,X^{t,x}_r,\tilde Y^{t,x,n}_r,\tilde Z^{t,x,n}_r)dr-n\int_t^T\!\!\int_E(\tilde V^{t,x,n}_r(e))^-\lambda(de)dr
    \\
    &\quad -\int_s^T \tilde Z^{t,x,n}_r dW_r-\int_{s}^T\!\!\!\int_E \tilde V^{t,x,n}_{r}(e)\mu(dr,de),\quad \forall s\in [t,T].
\end{align*}
By standard arguments, we find that
\begin{align*}
  \|Y^{t,x,0,n}\|^2_{\mcS^2_t}\leq C\E^{\nu^n}\Big[|\psi(X^{t,x}_T)|^2+\int_t^T |f(r,X^{t,x}_r,0,0)|^2dr\Big],
\end{align*}
and the result follows by Proposition~\ref{prop:SDEmoment}.\qed\\

%%%%%%%%%%%%%%%%%%%%%%%%%%%%%%%%%%%%%%%%%%%%%%%%%%%%%%%%%%%%%%%%%%%%%%%%%%%%%%%%%%%%%%%%%%%%%%%%%%%%%%%%%%%%%%%%%%%%%%%%%%%%%%%%%%%%%%%%%%%%%%%%%%
%%%%%%%%%%%%%%%%%%%%%%%%%%%%%%%%%%%%%%%%%%%%%%%%%%%%%%%%%%%%%%%%%%%%%%%%%%%%%%%%%%%%%%%%%%%%%%%%%%%%%%%%%%%%%%%%%%%%%%%%%%%%%%%%%%%%%%%%%%%%%%%%%%
%%%%%%%%%%%%%%%%%%%%%%%%%%%%%%%%%%%%%%%%%%%%%%%%%%%%%%%%%%%%%%%%%%%%%%%%%%%%%%%%%%%%%%%%%%%%%%%%%%%%%%%%%%%%%%%%%%%%%%%%%%%%%%%%%%%%%%%%%%%%%%%%%%

\section{Viscosity properties of the lower and upper game values\label{sec:visc}}

In this section, we begin by proving the first part of Theorem~\ref{thm:game-value}, demonstrating that the game introduced in the previous section has a value and that the corresponding value function is a viscosity solution to \eqref{ekv:var-ineq}. Our proof strategy hinges on a careful choice of candidate super- and sub-solutions, alongside a comparison principle for \eqref{ekv:var-ineq}, which is detailed in Appendix~\ref{app:uniq}. Specifically, we derive lower and upper bounds, $\underline v$ and $\bar v$, for the game value in \eqref{ekv:game-value}. We then show that $\underline v$ serves as a supersolution while $\bar v$ acts as a subsolution to \eqref{ekv:var-ineq}. Finally, using the comparison result presented in Proposition~\ref{app:prop:comp-visc}, we conclude that $\underline v=\bar v$ is the unique viscosity solution to \eqref{ekv:var-ineq}.

%In this section, we prove the first part of Theorem~\ref{thm:game-value}, showing that the game in the previous section has a value and that the value function is a viscosity solution to \eqref{ekv:var-ineq}. Our approach to prove this revolves around an intricate choice of candidate super- and sub-solutions along with a comparison principle for \eqref{ekv:var-ineq} that can be found in Appendix~\ref{app:uniq}. As candidate solutions we extract a lower and an upper bound for the game value in \eqref{ekv:game-value}. Below, we first derive the lower bound $\underline v$ and the upper bound $\bar v$ before we prove that $\underline v$ is a supersolution to \eqref{ekv:var-ineq} while $\bar v$ is a subsolution. A comparison result, presented in Proposition~\ref{app:prop:comp-visc}, then gives that $\underline v=\bar v$ is the unique viscosity solution to \eqref{ekv:var-ineq}.

\subsection{The lower bound}

%\fbox{
%\parbox{0.95\textwidth}{
%\begin{itemize}
%  \item For each $(t,x)\in [0,T]\times \R$ and $k\in\bbN$, the sequence $(Y^{t,x,k,n})_{n\in\bbN}$ is a non-increasing sequence of \cadlag processes. As the sequence is bounded from below by $\underline \mcY^{t,x}\in\mcS^2$ it follows that there is a process $\underline Y^{t,x,k}$ such that $Y^{t,x,k,n}\searrow \underline Y^{t,x,k}$ pointwisely.
%  \item Similarly, there is a map $\underline v_k:[0,T]\times \R^d\to\R$ such that $v_{k,n}\searrow \underline v_k$ and by pointwise convergence it holds that $\underline v_k(\eta,X^{t,x}_\eta)=\underline Y^{t,x,k}_\eta$ for each $\eta\in\mcT_t$.
%  \item We show that $\underline Y^{t,x,k}=\esssup_{\tau\in\mcT_s}Y^{t,x,k;\tau}_s$, where $(Y^{t,x,k;\tau},Z^{t,x,k;\tau},V^{t,x;\tau},K^{-,t,x,k;\tau})\in\mcS^{2,\tau} \times \mcH^{2,\tau}(W) \times \mcH^{2,\tau}(\mu) \times \mcA^{2,\tau}$ is the unique maximal solution to the BSDE
%\begin{align}\label{ekv:stopped-bsde}
%  \begin{cases}
%     Y^{t,x,k;\tau}_s=\supOP\underline v_{k-1}(\tau,X^{t,x}_\tau)+\int_s^\tau f(r,X^{t,x}_r, Y^{t,x,k;\tau}_r, Z^{t,x,k;\tau}_r)dr -\int_s^\tau  Z^{t,x,k;\tau}_r dW_r
%    \\
%    \quad-\int_{s}^\tau\!\!\!\int_E  V^{t,x,k;\tau}_{r}(e)\mu(dr,de)-(K^{-,t,x,k;\tau}_\tau- K^{-,t,x,k;\tau}_s),\quad \forall s\in [t,\tau],
%    \\
%    V^{t,x,k;\tau}_s(e)\geq - \chi(s,X^{t,x}_{s-},e),\quad d\Prob\otimes ds\otimes \lambda(de)-a.e.
%  \end{cases}
%\end{align}
%\end{itemize}
%}
%}

We introduce  the quadruple $(\underline Y^{t,x,0},\underline Z^{t,x,0},\underline V^{t,x,0},\underline K^{-,t,x,0}) \in\mcS^2_{t}\times \mcH^2_{t}(W) \times \mcH^2_{t}(\mu) \times \mcA^2_{t}$ defined as the unique maximal solution to the BSDE with constrained jumps,
\begin{align}\label{ekv:stopped-bsde-0}
  \begin{cases}
     \underline Y^{t,x,0}_s=\psi(X^{t,x}_T)+\int_s^T f(r,X^{t,x}_r, \underline Y^{t,x,0}_r, \underline Z^{t,x,0}_r)dr -\int_s^T  \underline Z^{t,x,0}_r dW_r-\int_{s}^T \!\!\!\int_E  \underline V^{t,x,0}_{r}(e)\mu(dr,de)
    \\
    \quad-(\underline K^{-,t,x,0}_T- \underline K^{-,t,x,0}_s),\quad \forall s\in [t,T],
    \\
    \underline V^{t,x,0}_s(e)\geq - \chi(s,X^{t,x}_{s-},e),\quad d\Prob\otimes ds\otimes \lambda(de)-a.e.,
  \end{cases}
\end{align}
and let $\underline v_0(t,x):=\underline Y^{t,x,0}_t$. Then combining Proposition~\ref{prop:BSDEmoment} with the results of \cite{Kharroubi2010} we find that $\underline v_0\in\Pi^g_c$. By repeated application of Theorem~\ref{thm:qvi-stop}, we can thus define a sequence $(\underline Y^{t,x,k})_{k\in\bbN}$ in $\mcS^2_t$ as $\underline Y^{t,x,k}_s=\esssup_{\tau\in\mcT_s}\underline Y^{t,x,k;\tau}_s$ for $k\geq 1$, where $(\underline Y^{t,x,k;\tau},\underline Z^{t,x,k;\tau},\underline V^{t,x,k;\tau},\underline K^{-,t,x,k;\tau}) \in\mcS^2_{[0,\tau]}\times \mcH^2_{[0,\tau]}(W) \times \mcH^2_{[0,\tau]}(\mu) \times \mcA^2_{[0,\tau]}$ is the unique maximal solution to the BSDE
\begin{align}\label{ekv:stopped-bsde-k}
  \begin{cases}
     \underline Y^{t,x,k;\tau}_s=\supOP \underline v_{k-1}(\tau,X^{t,x}_\tau)+\int_s^\tau f(r,X^{t,x}_r, \underline Y^{t,x,k;\tau}_r, \underline Z^{t,x,k;\tau}_r)dr -\int_s^\tau  \underline Z^{t,x,k;\tau}_r dW_r
    \\
    \quad-\int_{s}^\tau\!\!\!\int_E  \underline V^{t,x,k;\tau}_{r}(e)\mu(dr,de)-(\underline K^{-,t,x,k;\tau}_\tau- \underline K^{-,t,x,k;\tau}_s),\quad \forall s\in [0,\tau],
    \\
    \underline V^{t,x,k;\tau}_s(e)\geq - \chi(s,X^{t,x}_{s-},e),\quad d\Prob\otimes ds\otimes \lambda(de)-a.e.,
  \end{cases}
\end{align}
with $\underline v_{k-1}(t,x):=\underline Y^{t,x,k-1}_t$.

Furthermore, for each $(t,x)\in [0,T]\times\R^d$ and $u\in\bar\mcU_t$, we let $(\underline P^{t,x;u},\underline Q^{t,x;u},\underline R^{t,x;u},\underline S^{t,x;u})$ with $\underline P^{t,x;u}-\Xi^{t,x;u} \in \mcS^2$ and $(\underline Q^{t,x;u},\underline R^{t,x;u},\underline S^{t,x;u})\in\mcH^2(W)\times\mcH^2(\mu)\times \mcA^2$ be the maximal solution to the non-standard BSDE
\begin{align}
\begin{cases}
\underline P^{t,x;u}_s=\psi(X^{t,x;u}_T)+\int_s^T f(r,X^{t,x;u}_r,\underline P^{t,x;u}_r,\underline Q^{t,x;u}_r)dr-\int_s^T \underline Q^{t,x;u}_r dW_r - \int_{s}^T\!\!\!\int_E \underline R^{t,x;u}_{r}(e)\mu(dr,de)
\\
\quad - \Xi^{t,x;u}_{T}+\Xi^{t,x;u}_s-(\underline S^{t,x;u}_T- \underline S^{t,x;u}_s),\quad \forall s\in [t,T],
\\
\underline R^{t,x;u}_{s}(e)\geq - \chi(s,X^{t,x;u}_{s-},e),\quad d\Prob\otimes ds\otimes \lambda(de)-a.e.,
\end{cases}\label{ekv:non-ref-bsde-cstr-jmp}
\end{align}
the existence of which follows by repeating the argumentation in Section 3 of \cite{Kharroubi2010}. We also have that
\begin{align}\label{ekv:underP}
  \underline P^{t,x;u}_s=\essinf_{\nu\in\bar\mcV} P^{t,x;u,\nu}_s.
\end{align}

\begin{prop}\label{prop:dblunder}
For each $(t,x)\in [0,T]\times\R^d$, $k\in\bbN$ and $s\in [t,T]$, there is a $\underline u^{*,k}\in\mcU^{t,x,k}_s$ such that
\begin{align*}
  \underline Y^{t,x,k}_s\leq \underline P^{t,x;\underline u^{*,k}}_s.
\end{align*}
Moreover, $\underline v_k\in\Pi^g_c$ is the unique viscosity solution to
\begin{align}\label{ekv:dblunder-var-ineq}
\begin{cases}
  \min\{\underline v_k(t,x)-\supOP \underline v_{k-1}(t,x),\max\{\underline v_k(t,x)-\infOP \underline v_k(t,x),-\frac{\partial}{\partial t}\underline v_k(t,x)-\mcL \underline v_k(t,x)\\
  \quad-f(t,x,\underline v_k(t,x),\sigma^\top(t,x)\nabla_x \underline v_k(t,x))\}\}=0,\quad\forall (t,x)\in[0,T)\times \R^d \\
  \underline v_k(T,x)=\psi(x).
\end{cases}
\end{align}
\end{prop}

\noindent\emph{Proof.} That $(\underline v_k)_{k\in\bbN}$ is the unique viscosity solution to \eqref{ekv:dblunder-var-ineq} is immediate from the definition and Theorem~\ref{thm:qvi-stop}. To show the first part of the proposition, we let $\underline u^{*,k}:=(\tau^*_j,\beta^*_j)_{j=1}^{N^*}\in\mcU^{t,x,k}_s$ be defined through
\begin{itemize}
  \item $\tau^*_{j}:=\inf \big\{r \geq \tau^*_{j-1}:\:\underline v_{k-j+1}(r,X^{t,x;[\underline u^{*,k}]_{j-1}}_{r})=\supOP \underline v_{k-j}(r,X^{t,x;[\underline u^{*,k}]_{j-1}}_r)\big\}\wedge T$,
  \item $\beta^*_j\in\mathop{\arg\max}_{b\in U}\{\underline v_{k-j}(\tau^*_{j},X^{t,x;[\underline u^{*,k}]_{j-1}}_{\tau^*_{j}}+\xi(\tau^*_{j},X^{t,x;[\underline u^{*,k}]_{j-1}}_{\tau^*_{j}},b))-\ell(\tau^*_{j},X^{t,x;[\underline u^{*,k}]_{j-1}}_{\tau^*_{j}},b)\}$,
\end{itemize}
for $j=1,\ldots,k$ with $\tau_0^*:=s$, while $N^*:=\max\{j\in \{0,\ldots,k\}:\tau^*_j<T\}$.\\

We show that $\underline Y^{t,x,k}_s\leq P^{t,x;\underline u^{*,k}}_s$. By the definition of the pair $(\tau^*_1,\beta^*_1)$ we find that
\begin{align*}
  \underline Y^{t,x,k;\tau^*_1}_{s'}&= \supOP \underline v_{k-1}(\tau^*_1,X^{t,x}_{\tau_1^*}) + \int_{s'}^{\tau_1^*} f(r,X^{t,x}_r,\underline Y^{t,x,k;\tau^*_1}_r,\underline Z^{t,x,k;\tau^*_1}_r)dr
  \\
  &\quad - \int_{s'}^{\tau_1^*}Z^{t,x,k;\tau^*_1}_r dW_r-\int_{{s'}}^{\tau_1^*}\!\!\!\int_E \underline V^{t,x,k;\tau^*_1}_{r}(e)\mu(dr,de)-(\underline K^{-,t,x,k;\tau^*_1}_{\tau^*_1} - \underline K^{-,t,x,k;\tau^*_1}_{{s'}})
  \\
  &= \ett_{[\tau^*_1<T]}\{\underline v_{k-1}(\tau^*_1,X^{t,x;[\underline u^{*,k}]_1}_{\tau_1^*})-\ell({\tau_1^*},X^{t,x}_{\tau_1^*},\beta_1^*)\} +\ett_{[\tau^*_1=T]}\psi(X^{t,x}_T)
  \\
  &\quad +\int_{s'}^{\tau_1^*} f(r,X^{t,x}_r,\underline Y^{t,x,k;\tau^*_1}_r,\underline Z^{t,x,k;\tau^*_1}_r)dr - \int_{s'}^{\tau_1^*}\underline Z^{t,x,k;\tau^*_1}_r dW_r
  \\
  &\quad-\int_{{s'}}^{\tau_1^*}\!\!\!\int_E V^{t,x,k;\tau^*_1}_{r}(e)\mu(dr,de)-(\underline K^{-,t,x,k;\tau^*_1}_{\tau^*_1} - \underline K^{-,t,x,k;\tau^*_1}_{{s'}})
\end{align*}
for all ${s'}\in [t,\tau^*_1]$ and that $\underline V^{t,x,k;\tau^*_1}_{{s'}}(e)\geq - \chi({s'},X^{t,x}_{{s'}-},e)$, $d\Prob\otimes ds\otimes \lambda(de)-a.e.$. To simplify notation later, we let $(\mcY^0,\mcZ^0,\mcV^0,\mcK^0):=(\underline Y^{t,x,k;\tau^*_1},\underline Z^{t,x,k;\tau^*_1},\underline V^{t,x,k;\tau^*_1},\underline K^{-,t,x,k;\tau^*_1})$. Similarly, there is a quadruple $(\mcY^1,\mcZ^1,\mcV^1,\mcK^1)\in\mcS^2_t\times\mcH^2_t(W)\times\mcH^2_t(\mu)\times\mcA^2_{t}$ that is the maximal solution to the BSDE with constrained jumps
\begin{align*}
\begin{cases}
  \mcY^{1}_{s'}=\supOP \underline v_{k-2}(\tau^*_2,X^{t,x;[\underline u^{*,k}]_1}_{\tau_2^*})+\int_{s'}^{\tau_2^*} f\big(r,X^{t,x;[\underline u^{*,k}]_1}_r,\mcY^1_r,\mcZ^1_r\big)dr
  \\
  \quad-\int_{s'}^{\tau_2^*} \mcZ^1_r dW_r-\int_{{s'}}^{\tau_2^*}\!\!\!\int_E \mcV^{1}_{r}(e)\mu(dr,de)-(\mcK^1_{\tau_2^*}-\mcK^1_{s'}),\quad \forall {s'}\in [t,{\tau_2^*}],
\\
\mcV^{1}_{{s'}}(e)\geq - \chi({s'},X^{t,x;\underline u^{*,k}}_{{s'}-},e),\quad d\Prob\otimes ds'\otimes \lambda(de)-a.e.,
\end{cases}
\end{align*}
Arguing as in Section 3 of \cite{Kharroubi2010} (see also the proof of Theorem 3.1 in \cite{imp-stop-game}) we find that $\mcY^{1}=\lim_{n\to\infty}\mcY^{1,n}$, where $(\mcY^{1,n},\mcZ^{1,n},\mcV^{1,n})\in \mcS^2\times \mcH^2(W)\times\mcH^2(\mu)$ is the unique solution to the penalized stopped standard BSDE
\begin{align*}
  \mcY^{1,n}_{s'}&=\supOP \underline v_{k-2}(\tau^*_2,X^{t,x;[\underline u^{*,k}]_1}_{\tau_2^*})+\int_{s'}^{\tau_2^*} f^n\big(r,X^{t,x;[\underline u^{*,k}]_1}_r,\mcY^1_r,\mcZ^1_r,\mcV^1_r\big)dr-\int_{s'}^{\tau_2^*} \mcZ^1_r dW_r
  \\
  & \quad -\int_{{s'}}^{\tau_2^*}\!\!\!\int_E \mcV^{1}_{r}(e)\mu(dr,de)- (\mcK^1_{\tau_2^*}-\mcK^1_{s'}),\quad\forall {s'}\in[t,\tau^*_2]
\end{align*}
and continuity together with an approximation of $X^{t,x;[\underline u^{*,k}]_1}$ and the stability result for BSDEs with jumps in Proposition 2.2 of \cite{Barles97} implies that $\mcY^{1}_{{s'}}=\underline v_{k-1}({s'},X^{t,x;[\underline u^{*,k}]_1}_{{s'}})$, $\Prob$-a.s., for each ${s'}\in [\tau^*_1,T]$.

Repeating this process $k$ times we find that there is a sequence $(\mcY^j,\mcZ^j,\mcV^j,\mcK^j)_{j=0}^k\subset\mcS^2\times\mcH^2(W)\times\mcH^2(\mu)\times \mcA^2$ such that $\mcY:=\ett_{[t,\tau^*_{1}]}\mcY^0+\sum_{j=1}^k\ett_{(\tau^*_{j},\tau^*_{j+1}]}\mcY^j$, $\mcZ:=\sum_{j=0}^k\ett_{(\tau^*_{j},\tau^*_{j+1}]}\mcZ^j$,  $\mcV:=\sum_{j=0}^k\ett_{(\tau^*_{j},\tau^*_{j+1}]}\mcV^j$ and $\mcK_{s'}:=\sum_{j=0}^{k}\ett_{[\hat \tau_j<{s'}]}\{\mcK^j_{{s'}\wedge \tau^*_{j+1}}-\mcK^j_{\tau^*_j}\}$ satisfies
\begin{align*}
  \mcY_{s'}&=\psi(X^{t,x;\underline u^{*,k}}_T)+\int_{s'}^{T} f\big(r,X^{t,x;\underline u^{*,k}}_r,\mcY_r,\mcZ_r\big)dr-\int_{s'}^{T}\mcZ_rdW_r - \int_{{s'}}^{T}\!\!\!\int_E \mcV_{r}(e)\mu(dr,de)
  \\
  &\quad -\sum_{j=1}^{N^*}\ett_{[{s'}\leq \tau^*_j]}\ell({\tau_j^*},X^{t,x;[\underline u^{*,k}]_{j-1}}_{\tau_j^*},\beta_j^*)-(\mcK_T-\mcK_{s'}),\quad\forall {s'}\in[t,T].
\end{align*}
Since $P^{t,x;\underline u^{*,k}}$ is the maximal solution to this BSDE, we conclude that $\underline Y^{t,x,k}_s\leq P^{t,x;\underline u^{*,k}}_s$.\qed\\

Combining the above proposition with \eqref{ekv:underP} we find that:% and Proposition~\ref{prop:trunk-game} is the following:
%that $\underline Y^{t,x,k}\leq Y^{t,x,k,n}$ for each $k,n\in\bbN$. The following result implies that $\underline Y^{t,x,k}$ is also a lower bound for the value functions in Theorem~\ref{thm:game-value}
\begin{cor}
For each $(t,x)\in [0,T]\times\R^d$, $s\in[t,T]$ and $k\in\bbN$, it holds that
  \begin{align*}
    \underline Y^{t,x,k}_s\leq \esssup_{u\in\mcU^{t,x}_s}\essinf_{\nu\in\mcV^{t,x}}P^{t,x;u,\nu}_s.
  \end{align*}
\end{cor}

%\noindent\emph{Proof.} Since all members of $\mcV^S$ are uniformly bounded, we have, whenever $\nu^S\in\mcV^S$, that $\nu^S\in\mcV^{S,n}$ for $n$ sufficiently large. This implies that
%\begin{align*}
%  \underline Y^{t,x,k}_s\leq \underline Y^{t,x,k}_s=\essinf_{\nu^S\in\mcV^{S}_t}\esssup_{u\in\mcU^{t,k}_s}P^{t,x;u,\nu^S(u)}_s\leq \essinf_{\nu^S\in\mcV^{S}_t}\esssup_{u\in\mcU^t_s}P^{t,x;u,\nu^S(u)}_s.
%\end{align*}
%On the other hand, it follows immediately from Proposition~\ref{prop:dblunder} that
%  \begin{align*}
%    \underline Y^{t,x,k}_s=  \essinf_{\nu\in\mcV}P^{t,x;u^{*,k},\nu}_s \leq \esssup_{u\in\mcU_s}\essinf_{\nu\in\mcV}P^{t,x;u,\nu}_s,
%  \end{align*}
%which concludes the proof.\qed\\

A comparison result for the corresponding BSDEs gives that the sequence $(\underline v_k)_{k\in\bbN}$ is non-decreasing. We thus conclude that there is a map $\underline v:[0,T]\times\R^d\to\R$ such that $\underline v_k\nearrow \underline v$ pointwisely and since each $\underline v_k$ is continuous it follows that $\underline v$ is l.s.c.

\subsection{The upper bound}

%\fbox{
%\parbox{0.85\textwidth}{
%\begin{itemize}
%  \item For each $(t,x)\in [0,T]\times \R$ and $n\in\bbN$, the sequence $(Y^{t,x,k,n})_{k\in\bbN}$ is a non-decreasing sequence of \cadlag processes. As the sequence is bounded from above by $\bar \mcY^{t,x}\in\mcS^2$ it follows that there is a process $\bar Y^{t,x,n}$ such that $Y^{t,x,k,n}\nearrow \bar Y^{t,x,n}$ pointwisely.
%  \item Similarly, there is a map $\bar v_n:[0,T]\times \R^d\to\R$ such that $v_{k,n}\nearrow \bar v_n$ and by pointwise convergence it holds that $\bar v_n(\eta,X^{t,x}_\eta)=\bar Y^{t,x,n}_\eta$ for each $\eta\in\mcT_t$.
%  \item $\underline v_n$ is the unique viscosity solution to
%  \begin{align}\label{ekv:var-ineq-n}
%\begin{cases}
%  \min\{v_{n}(t,x)-\supOP v_{n}(t,x),-\frac{\partial}{\partial t}v_{n}(t,x)-\mcL v_{n}(t,x)-\mcK^n v_{n}(t,x)\\
%\quad-f(t,x,v_{n}(t,x),\sigma^\top(t,x)\nabla_x v_{n}(t,x))\}=0,  \quad\forall (t,x)\in[0,T)\times \R^d \\
%  v_{n}(T,x)=\psi(x).
%\end{cases}
%\end{align}
%\end{itemize}
%These results should follow from \cite{qvi-Levy}.
%}
%}

For each $n\in\bbN$, let $\bar v_n\in\Pi^g_c$ be the unique viscosity solution in $\Pi^g_c$ to \eqref{ekv:var-ineq-n} and let the quadruple $(\bar Y^{t,x,n},\bar Z^{t,x,n},\bar V^{t,x,n},\bar K^{+,t,x,n})\in\mcS^{2}_t \times \mcH^{2}_t(W) \times \mcH^{2}_t(\mu) \times \mcA^{2}_t$ be the unique solution to the reflected BSDE
\begin{align}\label{ekv:rbsde-jmp}
  \begin{cases}
     \bar Y^{t,x,n}_s=\psi(X^{t,x,n}_T)+\int_s^T f^n(r,X^{t,x,n}_r, \bar Y^{t,x,n}_r, \bar Z^{t,x,n}_r, \bar V^{t,x,n}_{r})dr -\int_s^T  \bar Z^{t,x,n}_r dW_r-\int_{s}^T\!\!\!\int_E \bar V^{t,x,n}_{r}(e)\mu(dr,de)
    \\
    \quad+\bar K^{+,t,x,n}_T - \bar K^{+,t,x,n}_s,\quad \forall s\in [t,T],
    \\
    \bar Y^{t,x,n}_s\geq \supOP \bar v_n(s,X^{t,x,n}_s),\, \forall s\in [t,T] \quad\text{and}\quad\int_t^T\!\! \big(\bar Y^{t,x,n}_s-\supOP \bar v_n(s,X^{t,x,n}_s)\big)d\bar K^{+,t,x,n}_s.
  \end{cases}
\end{align}
By Proposition~\ref{prop:rbsde-jmp}, $\tilde v\in\Pi^g_c$ defined as $\tilde v(t,x)=\bar Y^{t,x,n}_t$ for each $(t,x)\in[0,T]\times\R^d$ is the unique viscosity solution in $\Pi^g_c$ to
\begin{align}\label{ekv:qvi-altered}
\begin{cases}
  \min\{\tilde v(t,x)-\supOP \bar v_n(t,x),-\tilde v_t(t,x)-\mcL \tilde v(t,x)-\mcK^n \tilde v(t,x)\\
  \quad-f(t,x,\tilde v(t,x),\sigma^\top(t,x)\nabla_x \tilde v(t,x))\}=0,\quad\forall (t,x)\in[0,T)\times \R^d \\
  \tilde v(T,x)=\psi(x).
\end{cases}
\end{align}
On the other hand, $\bar v_n$ is also a solution to \eqref{ekv:qvi-altered} and we conclude that $\tilde v=\bar v_n$. In particular, we find that $\bar Y^{t,x,n}_\eta=\bar v_n(\eta,X^{t,x}_\eta)$ for any stopping time $\eta\in\mcT_t$. We have the following alternative representation of $\bar Y^{t,x,n}$:

\begin{prop}\label{prop:trunk-game-value-upper}
For each $(t,x)\in[0,T]\times\R^d$, $s\in[t,T]$ and $n\in\bbN$, we have
\begin{align}\label{ekv:trunk-game-value-upper}
  \bar Y^{t,x,n}_s \geq \essinf_{\nu\in\mcV^{t,x}}\esssup_{u\in\mcU^{t,x}_s} P^{t,x;u,\nu}_s.%=\esssup_{u\in\mcU^{t,x}_s}\essinf_{\nu\in\mcV^{t,x,n}}P^{t,x;u,\nu}_s
\end{align}
\end{prop}

\noindent\emph{Proof.} We let $\nu^{n,*}\in\mcV^{t,x,n}$ be defined as
\begin{align*}
  &\nu^{n,*}(u)(r,e):=n\ett_{\bar A^{t,x,u}_r(e)}
\end{align*}
for all $(u,r,e)\in \bar\mcU_s\times [0,T]\times E$, where
\begin{align*}
  \bar A^{t,x,u}_r(e):=&\{(r,\omega,e)\in [0,T]\times \Omega \times E: \bar v_n(r,X^{t,x;u}_{r-}+\xi(r,X^{t,x;u}_{r-},e))+\chi(r,X^{t,x;u}_{r-},e)
  \\
  &\quad-\bar v_n(r,X^{t,x;u}_{r-} + \xi(r,X^{t,x;u}_{r-},e))<0\}.
\end{align*}
Remark~\ref{rem:Vtx-cont} then implies that
\begin{align*}
  \nu^{n,*}(u)(r,e)=n\ett_{[0,\tau_1]}(r)\ett_{[\bar V^{t,x,n}_{r}(e)<-\chi(r,X^{t,x}_{r},e)]}+n\sum_{j=1}^\infty\ett_{(\tau_j,\tau_{j+1}]}(r) \ett_{[\bar V^{\tau_j,X^{t,x;u}_{\tau_j},n}_{r}(e)<-\chi(r,X^{t,x;u}_{r},e)]}
\end{align*}
$d\Prob\otimes\lambda(de)\otimes dr$-a.e. Let $u^*:=(\tau^*_j,\beta^*_j)_{j=1}^{N^*}$ be defined through
\begin{itemize}
  \item $\tau^*_{j}:=\inf \big\{{s'} \geq \tau^*_{j-1}:\:v({s'},X^{t,x;[u^*]_{j-1}}_{{s'}})=\supOP v({s'},X^{t,x;[u^*]_{j-1}}_{s'})\big\}\wedge T$,
  \item $\beta^*_j\in\mathop{\arg\max}_{b\in U}\{v(\tau^*_{j},X^{t,x;[u^*]_{j-1}}_{\tau^*_{j}}+\xi(\tau^*_{j},X^{t,x;[u^*]_{j-1}}_{\tau^*_{j}},b)) - \ell(\tau^*_{j},X^{t,x;[u^*]_{j-1}}_{\tau^*_{j}},b)\}$,
\end{itemize}
for $j=1,\ldots$ with $\tau_0^*:=s$, while $N^*:=\max\{j\in \bbN:\tau^*_j<T\}$. In the proof, which is divided into three steps, we show that $u^*\in\mcU^{t,x}_s$ and that $\bar Y^{t,x,n}_s=P^{t,x;u^*,\nu^{n,*}(u^*)}_s$, where in the latter $u^*$ is an optimal response to $\nu^{n,*}$.\\

\textbf{Step 1:} We first show that there exists of a triple $(\mcY,\mcZ,\mcV)$ such that for each $k\in\bbN$, we have that $(\mcY,\mcZ,\mcV)\in\mcS^{2}_{[t,\tau^*_{k+1}]} \times \mcH^{2}_{[t,\tau^*_{k+1}]}(W) \times \mcH^{2}_{[t,\tau^*_{k+1}]}(\mu)$ satisfies
%\begin{align*}
%  \mcY_{s'}=\ett_{[t,\tau^*_1]}({s'})\bar Y^{t,x,n}_{{s'}} + \sum_{j=1}^k \ett_{(\tau^*_j,\tau^*_{j+1}]}({s'})\bar Y^{\tau^*_j,X^{t,x;[u^*]_j}_{\tau_j},n}_{s'}=\bar v_n({s'},X^{t,x;u^*_{{s'}-}}_{s'}),\quad\forall {s'}\in[t,\tau^*_{k+1}]
%\end{align*}
%and
\begin{align}\nonumber
  \mcY_{s'}&=\supOP \bar v_n(\tau^*_{k+1},X^{t,x;[u^*]_{k}}_{\tau^*_{k+1}})+\int_{s'}^{\tau^*_{k+1}} f^n\big(r,X^{t,x;u^*}_r,\mcY_r,\mcZ_r,\mcV_r\big)dr-\int_{s'}^{\tau^*_{k+1}}\mcZ_rdW_r
  \\
  &\quad - \int_{{s'}}^{\tau^*_{k+1}}\!\!\!\int_E \mcV_{r}(e)\mu(dr,de) -\sum_{j=1}^{N^*\wedge k}\ett_{[{s'}\leq \tau^*_j]}\ell({\tau_j^*},X^{t,x;[u^*]_{j-1}}_{\tau_j^*},\beta_j^*),\quad\forall {s'}\in[t,\tau^*_{k+1}].\label{ekv:bsde-k}
\end{align}
Actually, since $f^n$ is Lipschitz in $(y,z,v)$, uniformly in $(t,x)$, standard results for BSDEs with jumps implies that if \eqref{ekv:bsde-k} holds, then there is a constant $C>0$ such that
\begin{align}\nonumber
  &\|\mcY\|^2_{\mcS^2_{[t,\tau^*_{k+1}]}}+\|\mcZ\|^2_{\mcH^2_{[t,\tau^*_{k+1}]}(W)}+\|\mcV\|^2_{\mcH^2_{[t,\tau^*_{k+1}]}(\mu)}
  \\
  &\leq C\E\Big[\supOP \bar v^2_n(\tau^*_{k+1},X^{t,x;[u^*]_{k}}_{\tau^*_{k+1}})+\int_t^{\tau^*_{k+1}} |f^n\big(r,X^{t,x;u^*}_r,0,0,0\big)|^2dr
  +\Big(\sum_{j=1}^{N^*\wedge k}\ell({\tau_j^*},X^{t,x;[u^*]_{j-1}}_{\tau_j^*},\beta_j^*)\Big)^2\Big]\nonumber
  \\
  &\leq C\Big(1+|x|^\rho+\E\Big[\Big(\sum_{j=1}^{N^*}\ell({\tau_j^*},X^{t,x;[u^*]_{j-1}}_{\tau_j^*},\beta_j^*)\Big)^2\Big]\Big),\label{ekv:yzv-bnd-k}
\end{align}
for all $(t,x,k)\in[0,T]\times\R^d\times\bbN$. Hence, we only need to show existence of a suitably measurable triple $(\mcY,\mcZ,\mcV)$ satisfying \eqref{ekv:bsde-k}.

To begin with, we observe that $\bar K^{+,t,x,n}_{\tau^*_1}-\bar K^{+,t,x,n}_{t}=0$, $\Prob$-a.s., and we find that
\begin{align*}
  \bar Y^{t,x,n}_{s'}&= \supOP \bar v_{n}(\tau^*_1,X^{t,x}_{\tau_1^*}) + \int_{s'}^{\tau_1^*} f^n(r,X^{t,x}_r,\bar Y^{t,x,n}_r,\bar Z^{t,x,n}_r,\bar V^{t,x,n}_{r})dr
  \\
  &\quad - \int_{s'}^{\tau_1^*}\bar Z^{t,x,n}_r dW_r-\int_{{s'}}^{\tau_1^*}\!\!\!\int_E \bar V^{t,x,n}_{r}(e)\mu(dr,de)
  \\
  &= \ett_{[\tau^*_1<T]}\{\bar v_{n}(\tau^*_1,X^{t,x;[u^{*}]_1}_{\tau_1^*})-\ell({\tau_1^*},X^{t,x}_{\tau_1^*},\beta_1^*)\} +\ett_{[\tau^*_1=T]}\psi(X^{t,x}_T)
  \\
  &\quad +\int_{s'}^{\tau_1^*} f^n(r,X^{t,x}_r,\bar Y^{t,x,n}_r,\bar Z^{t,x,n}_r,\bar V^{t,x,n}_{r})dr - \int_{s'}^{\tau_1^*}\bar Z^{t,x,n}_r dW_r -\int_{{s'}}^{\tau_1^*}\!\!\!\int_E \bar V^{t,x,n}_{r}(e)\mu(dr,de)
\end{align*}
for all ${s'}\in [t,\tau^*_1]$. To simplify notation later, we let $(\mcY^0,\mcZ^0,\mcV^0):=(\bar Y^{t,x,n},Z^{t,x,n},\bar V^{t,x,n})$. Similarly, there is a quadruple $(\mcY^1,\mcZ^1,\mcV^1,\mcK^1)\in\mcS^2_t\times\mcH^2_t(W)\times\mcH^2_t(\mu)\times\mcA^2_{t}$ that solves the reflected bsde
\begin{align*}
\begin{cases}
  \mcY^{1}_{s'}=\psi(X^{t,x;[u^{*}]_1}_T)+\int_{s'}^T f^n\big(r,X^{t,x;[u^{*}]_1}_r,\mcY^1_r,\mcZ^1_r,\mcV^1_r\big)dr-\int_{s'}^T \mcZ^1_r dW_r-\int_{{s'}}^{T}\!\!\!\int_E \mcV^{1}_{r}(e)\mu(dr,de)
  \\
  \quad+ \mcK^1_T-\mcK^1_{s'},\quad\forall {s'}\in[\tau^*_1,T], \\
  \mcY^{1}_{s'}\geq \supOP v_{n}({s'},X^{t,x;[u^{*}]_1}_{s'}),\:\forall {s'}\in[\tau^*_1,T]\quad{\rm and}\quad
  \int_{\tau_1^*}^T(\mcY^{1}_{s'}-\supOP \bar v_{n}({s'},X^{t,x;[u^{*}]_1}_{s'}))d\mcK^{1}_{s'}=0.
\end{cases}
\end{align*}
Continuity together with an approximation of $X^{t,x;[u^{*}]_1}$ and a stability result for reflected BSDEs with jumps akin to the one in Proposition 3.6 of \cite{ElKaroui1} (see \eg Proposition A.1 in \cite{Dumitrescu15}) implies that $\mcY^{1}_{{s'}}=\bar v_{n}({s'},X^{t,x;[u^{*}]_1}_{{s'}})$, $\Prob$-a.s., for each ${s'}\in [\tau^*_1,T]$. In particular, as Proposition~\ref{prop:rbsde-jmp} gives that $\mcK^1_{\tau^*_2}-\mcK^1_{\tau^*_1}=0$, $\Prob$-a.s., we conclude that
\begin{align*}
  \mcY^{1}_{s'}&= \ett_{[\tau^*_2<T]}\{\bar v_{n}(\tau^*_2,X^{t,x;[u^{*}]_2}_{\tau_2^*})-\ell({\tau_2^*},X^{t,x;[u^{*}]_1}_{\tau_2^*},\beta_2^*)\} +\ett_{[\tau^*_2=T]}\psi(X^{t,x;u^{*}}_T)
  \\
  &\quad \int_{s'}^{\tau_2^*} f^n(r,X^{t,x}_r,\mcY^{1}_r,\mcZ^{1}_r,\mcV^{1}_r)dr - \int_{s'}^{\tau_2^*}\mcZ^{1}_r dW_r - \int_{{s'}}^{\tau_2^*}\!\!\!\int_E \mcV^{1}_{r}(e)\mu(dr,de)
\end{align*}
for all ${s'}\in [\tau^*_1,\tau^*_2]$.

Repeating this process $k$ times we find that there is a sequence $(\mcY^j,\mcZ^j,\mcV^j)_{j=0}^k\subset\mcS^2\times\mcH^2(W)\times\mcH^2(\mu)$ such that $\mcY:=\ett_{[t,\tau^*_{1}]}\mcY^0+\sum_{j=1}^k\ett_{(\tau^*_{j},\tau^*_{j+1}]}\mcY^j$, $\mcZ:=\sum_{j=0}^k\ett_{(\tau^*_{j},\tau^*_{j+1}]}\mcZ^j$ and $\mcV:=\sum_{j=0}^k\ett_{(\tau^*_{j},\tau^*_{j+1}]}\mcV^j$ satisfies \eqref{ekv:bsde-k}.\\
%\begin{align*}
%  \mcY_{s'}&=\psi(X^{t,x;u^{*}}_T)+\int_{s'}^{T} f^n\big(r,X^{t,x;u^{*}}_r,\mcY_r,\mcZ_r,\mcV_r\big)dr-\int_{s'}^{T}\mcZ_rdW_r - \int_{{s'}}^{T}\!\!\!\int_E \mcV_{r}(e)\mu(dr,de)
%  \\
%  &\quad -\sum_{j=1}^{N^*}\ett_{[{s'}\leq \tau^*_j]}\ell({\tau_j^*},X^{t,x;[u^{*}]_{j-1}}_{\tau_j^*},\beta_j^*),\quad\forall {s'}\in[t,T].
%\end{align*}
%Now, $\mcY-\Xi^{t,x;u}\in\mcS^2_t$ and $(\mcZ,\mcV)\in\mcH^2_t(W)\times\mcH^2_t(\mu)$ and by uniqueness of solutions to \eqref{ekv:non-ref-bsde-simp} we conclude that $Y^{t,x,k,n}_s=P^{t,x;u^{*},\nu^{n,*}}_s$.\\

\textbf{Step 2:} We prove that
\begin{align}\label{ekv:Xi-in-L2}
  \E\Big[\Big(\sum_{j=1}^{N^*}\ell({\tau_j^*},X^{t,x;[u^*]_{j-1}}_{\tau_j^*},\beta_j^*)\Big)^2\Big]<\infty
\end{align}
which, since $\ell\geq\delta>0$, implies that $N^*<\infty$, $\Prob$-a.s. and we can take the limit as $k\to\infty$ in \eqref{ekv:yzv-bnd-k} to conclude that $(\mcY,\mcZ,\mcV)\in\mcS^{2}_{t} \times \mcH^{2}_{t}(W) \times \mcH^{2}_{t}(\mu)$ whereas taking the limit as $k\to\infty$ in \eqref{ekv:bsde-k} and using that $\supOP \bar v_n(T,x)=\psi(T,x)$, gives
\begin{align*}
  \mcY_{s'}&=\psi(X^{t,x;u^*}_T)+\int_{s'}^T f^n\big(r,X^{t,x;u^*}_r,\mcY_r,\mcZ_r,\mcV_r\big)dr-\int_{s'}^T\mcZ_rdW_r - \int_{{s'}}^T\!\!\!\int_E \mcV_{r}(e)\mu(dr,de)
  \\
  &\quad -\sum_{j=1}^{N^*}\ett_{[{s'}\leq \tau^*_j]}\ell({\tau_j^*},X^{t,x;[u^*]_{j-1}}_{\tau_j^*},\beta_j^*),\quad\forall s'\in[t,T].
\end{align*}
Finally, \eqref{ekv:Xi-in-L2} guarantees that $u^*\in\mcU^{t,x}_s$ and by uniqueness of solutions to \eqref{ekv:non-ref-bsde-simp} it follows that  $\bar Y^{t,x,n}_s=P^{t,x;u^*,\nu^{n,*}(u^*)}_s$.

Fix $k\in\bbN$ and note that since $f^n(t,x,y,z,v)\leq f(t,x,y,z)$, comparison gives that $\mcY_{s'}\leq \tilde \mcY^k_{s'}$ for all $s'\in [t,\tau^*_{k+1}]$, where $(\tilde \mcY^k,\tilde \mcZ^k,\tilde \mcV^k)\in \mcS^{2}_{[t,\tau^*_{k+1}]} \times \mcH^{2}_{[t,\tau^*_{k+1}]}(W) \times \mcH^{2}_{[t,\tau^*_{k+1}]}(\mu)$ solves
\begin{align}\nonumber
  \tilde \mcY^k_{s'}&=\supOP \bar v_n(\tau^*_{k+1},X^{t,x;[u^*]_{k}}_{\tau^*_{k+1}})+\int_{s'}^{\tau^*_{k+1}} f\big(r,X^{t,x;u^*}_r,\tilde \mcY^k_r,\tilde \mcZ^k_r\big)dr-\int_{s'}^{\tau^*_{k+1}}\tilde \mcZ^k_rdW_r - \int_{{s'}}^{\tau^*_{k+1}}\!\!\!\int_E \tilde \mcV^k_{r}(e)\mu(dr,de)
  \\
  &\quad -\sum_{j=1}^{N^*\wedge k}\ett_{[{s'}\leq \tau^*_j]}\ell({\tau_j^*},X^{t,x;[u^*]_{j-1}}_{\tau_j^*},\beta_j^*),\quad\forall {s'}\in[t,\tau^*_{k+1}].\label{ekv:tilde-mcY}
\end{align}
To simplify notation we let $\bar Y_k:=\supOP \bar v_n(\tau^*_{k+1},X^{t,x;[u^*]_{k}}_{\tau^*_{k+1}})$, $X=X^{t,x;u^*}$ and $X^j=X^{t,x;[u^*]_j}$. Letting\footnote{Using the convention that $\frac{0}{0}0=0$.}
\begin{align*}
 \zeta^k_1(r):=\frac{f(r,X_r,\tilde \mcY^k_r,\tilde \mcZ^k_r) -f(r,X_r,0,\tilde \mcZ^k_r)}{\tilde \mcY^k_r}\ett_{[\tilde \mcY^k_r\neq 0]}
\end{align*}
and
\begin{align*}
 \zeta^k_2(r):=\frac{f(r,X_r,0,\tilde \mcZ^k_r) - f(r,X_r,0,0)}{|\tilde \mcZ^k_r|^2}(\tilde \mcZ^k_r)^\top
\end{align*}
we have by the Lipschitz continuity of $f$ that $|\zeta^k_1(r)|\vee|\zeta^k_2(r)|\leq k_f$. Using It\^o's formula, we find that
\begin{align*}
\tilde \mcY^k_{s'}&=R^k_{{s'},\tau^*_{k+1}}\bar Y_k+\int_{s'}^{\tau^*_{k+1}}R^k_{{s'},r}f(r,X_r,0,0)dr -\int_{s'}^{\tau^*_{k+1}} R^k_{{s'},r}\tilde\mcZ_rdW_r-\sum_{j=1}^{N^*\wedge k} R^k_{{s'},\tau^*_{j}}\ett_{[\tau^*_j\geq {s'}]}\ell(\tau^*_j,X^{j-1}_{\tau^*_{j}},\beta^*_j)
\end{align*}
with $R^k_{{s'},r}:=e^{\int_{s'}^{r}(\zeta^k_1(v)-\frac{1}{2}|\zeta^k_2(v)|^2)dv+\frac{1}{2}\int_{s'}^{r}\zeta^k_2(v)dW_v}$. Since the intervention costs are positive, taking the conditional expectation on both sides and using Proposition~\ref{prop:SDEmoment} gives
\begin{align*}
\tilde \mcY^k_{s'}&\leq\E\Big[R^k_{{s'},\tau^*_{k+1}}\bar Y_k+\int_{s'}^{\tau^*_{k+1}}R^k_{{s'},r}f(r,X_r,0,0)dr \Big|\mcF^t_{s'}\Big]
\\
&\leq C\Big(1+\E\big[(R^k_{{s'},\tau^*_{k+1}})^2\big|\mcF^t_{s'}\big]^{1/2}\E\Big[\sup_{r\in [{s'},\tau^*_{k+1}]}|X_r|^{2\rho} \Big|\mcF^t_{s'}\Big]^{1/2}\Big)
\\
&\leq C(1+|X_{s'}|^\rho),
\end{align*}
where the constant $C>0$ is independent of $(t,{s'},x)$ and $k$. On the other hand, there is another constant $C>0$ (that is also independent of $(t,{s'},x)$ and $k$) such that
\begin{align*}
  \tilde \mcY^k_{s'}\geq \mcY_{s'}\geq  \supOP\bar v^{n}({s'},X^{t,x;u^*_{{s'}-}}_{s'}) \geq -C(1+|X^{t,x;u^*_{{s'}-}}_{s'}|^\rho),
\end{align*}
$\Prob$-a.s., for all ${s'}\in [t,\tau^*_{k+1}]$. Proposition~\ref{prop:SDEmoment} then gives
\begin{align*}
\E\Big[\sup_{{s'}\in[t,\tau^*_{k+1}]}|\tilde \mcY^k_{s'}|^2\Big]&\leq C(1+|x|^{2\rho}),
\end{align*}
where $C>0$ does not depend on $(t,x)$ and $k$. Next, we derive a bound on the $\mcH^2_{[t,\tau^*_{k+1}]}(W)$ (resp. $\mcH^2_{[t,\tau^*_{k+1}]}(\mu)$) norm of $\tilde\mcZ^k$ (resp. $\tilde\mcV^k$) that is uniform in $k$. Applying It{\^o}'s formula to $|\tilde \mcY^k|^{2}$ we get
\begin{align}\nonumber
|\tilde \mcY^k_{\tau^*_{k+1}}|^2&= |\tilde \mcY^k_{t}|^2 - 2\int_{t}^{\tau^*_{k+1}} \tilde \mcY^k_r f(r,X_r,\tilde\mcY^k_r,\tilde\mcZ^k_r)dr + 2\int_{t}^{\tau^*_{k+1}}\tilde\mcY^k_r\tilde\mcZ^k_rdW_r + \int_{t}^{\tau^*_{k+1}}|\tilde\mcZ^k_r|^2dr
  \\
  &\quad+\int_{t}^{\tau^*_{k+1}}\!\!\!\int_E (2\tilde\mcY^k_{r-}\tilde\mcV^k_{r}(e)+|\tilde\mcV^k_{r}(e)|^2)\mu(dr,de) + \sum_{j=1}^{N^*}(2\tilde\mcY^{k,j-1}_{\tau^*_{j}}\ell(\tau^*_j,X^{j-1}_{\tau^*_{j}},\beta^*_j) + \ell^2(\tau^*_j,X^{j-1}_{\tau^*_{j}},\beta^*_j))
%|\mcY_{t}|^{2}+\int_{t}^T|\mcZ_s|^2ds+\int_{t}^T&=\psi^2(X_T)+2\int_t^T P_sf(s,X_s,P_s,Q_s)ds-2\int_t^T P_sQ_sdW_s
\label{ekv:from-ito}
\end{align}
where $\tilde\mcY^{k,j-1}$ is $\tilde\mcY^k$ without the $j-1$ first intervention costs. Since the intervention costs are nonnegative, we have
\begin{align*}
-\sum_{j=1}^{N^*} (2\tilde\mcY^{k,j-1}_{\tau^*_{j}}\ell(\tau^*_j,X^{j-1}_{\tau^*_{j}},\beta^*_j)+\ell^2(\tau^*_j,X^{j-1}_{\tau^*_{j}},\beta^*_j)) &\leq 2\sup_{{s'}\in [t,\tau^*_{k+1}]}|\tilde\mcY^k_{{s'}}|\sum_{j=1}^{N^*} \ell(\tau^*_j,X^{j-1}_{\tau^*_{j}},\beta^*_j)
\\
&\leq \kappa \sup_{{s'}\in [t,\tau^*_{k+1}]}|\tilde\mcY^k_{{s'}}|^2+\frac{1}{\kappa}\Big(\sum_{j=1}^{N^*\wedge k} \ell(\tau^*_j,X^{j-1}_{\tau^*_{j}},\beta^*_j)\Big)^2
\end{align*}
for any $\kappa>0$. Moreover,
\begin{align*}
  -2\int_{t}^{\tau^*_{k+1}}\!\!\!\int_E \tilde\mcY^k_{r-}\tilde\mcV^k_{r}(e)\mu(dr,de)\leq 2\int_{t}^{\tau^*_{k+1}}\!\!\!\int_E|\tilde\mcY^k_{r-}|^2\mu(dr,de)+\frac{1}{2}\int_{t}^{\tau^*_{k+1}}\!\!\!\int_E |\tilde\mcV^k_{r}(e)|^2\mu(dr,de)
\end{align*}
Inserted in \eqref{ekv:from-ito} and using the Lipschitz property of $f$ this gives
\begin{align}\nonumber
&|\tilde\mcY^k_{t}|^{2}+\frac{1}{2}\int_{t}^{\tau^*_{k+1}}| \tilde\mcZ^k_r|^2dr+\frac{1}{2}\int_{t}^{\tau^*_{k+1}}\!\!\!\int_E |\tilde\mcV^k_{r}(e)|^2\mu(dr,de)
\\
&\leq (C+\kappa)\sup_{s'\in[t,{\tau^*_{k+1}}]}|\tilde\mcY^k_{s'}|^2+\int_{t}^{\tau^*_{k+1}}|f(r,X_r,0,0)|^2dr + 2\int_{t}^{\tau^*_{k+1}}\!\!\!\int_E|\tilde\mcY^k_{r-}|^2\mu(dr,de)\nonumber
\\
&\quad-2\int_{t}^{\tau^*_{k+1}} \tilde\mcY^k_r\tilde\mcZ^k_rdW_r+\frac{1}{\kappa}\Big(\sum_{j=1}^{N^*\wedge k} \ell(\tau^*_j,X^{j-1}_{\tau^*_{j}},\beta^*_j)\Big)^2.\label{ekv:PQ-bound}
\end{align}
Now, as $[u^*]_k\in\mcU^k$, it follows that the stochastic integral is uniformly integrable and thus a martingale. To see this, note that the Burkholder-Davis-Gundy inequality gives
\begin{align*}
\E\Big[\sup_{{s'}\in[t,{\tau^*_{k+1}}]}\Big|\int_{t}^{s'} \tilde\mcY^k_r\tilde\mcZ^k_rdW_r\Big|\Big]\leq C\E\Big[\Big(\int_{t}^{\tau^*_{k+1}} |\tilde\mcY^k_{r}\tilde\mcZ^k_r|^2dr\Big)^{1/2}\Big]\leq C\E\Big[\sup_{s'\in [t,{\tau^*_{k+1}}]}|\tilde\mcY^k_{s'}|^2+\int_{t}^{\tau^*_{k+1}} |\tilde\mcZ^k_r|^2dr\Big]
\end{align*}
where the right-hand side is finite. Taking expectations on both sides of \eqref{ekv:PQ-bound} thus gives
\begin{align*}
\E\Big[\int_{t}^{\tau^*_{k+1}}| \tilde\mcZ^k_r|^2dr+\int_{t}^{\tau^*_{k+1}}\!\!\int_E| \tilde\mcV^k_r(e)|^2\lambda(de)dr\Big]&\leq C(1+\kappa)(1+|x|^{2\rho})+\frac{2}{\kappa}\E\Big[\Big(\sum_{j=1}^{N^*\wedge k} \ell(\tau^*_j,X^{j-1}_{\tau^*_{j}},\beta^*_j)\Big)^2\Big].
\end{align*}
Rearranging \eqref{ekv:tilde-mcY} and squaring both sides then gives
\begin{align*}
&\E\Big[\Big(\sum_{j=1}^{N^*\wedge k} \ell(\tau^*_j,X^{j-1}_{\tau^*_{j}},\beta^*_j)\Big)^2\Big]
\\
&\leq C\E\Big[|\tilde\mcY^k_{t}|^2+|\tilde\mcY^k_{\tau^*_{k+1}}|^2+\int_{t}^{\tau^*_{k+1}}|f(r,X_r,\tilde\mcY^k_r,\tilde\mcZ^k_r)|^2dr +\int_{t}^{\tau^*_{k+1}}|\tilde\mcZ^k_r|^2dr+\int_{t}^{\tau^*_{k+1}}\!\!\int_E| \tilde\mcV^k_r(e)|^2\lambda(de)dr\Big]
\\
&\leq C\E\Big[|\tilde\mcY_{\tau^*_{k+1}}|^2+\sup_{s'\in[{t},{\tau^*_{k+1}}]}|\tilde\mcY^k_{s'}|^2+\int_{t}^{\tau^*_{k+1}}(|f(r,X_r,0,0)|^2+|\tilde\mcZ^k_r|^2 + \int_E| \tilde\mcV^k_r(e)|^2\lambda(de))dr\Big]
\\
&\leq C(1+\kappa)(1+|x|^{2\rho})+\frac{C}{\kappa}\E\Big[\Big(\sum_{j=1}^{N^*\wedge k} \ell(\tau^*_j,X^{j-1}_{\tau^*_{j}},\beta^*_j)\Big)^2\Big].
\end{align*}
We can now prove that \eqref{ekv:Xi-in-L2} holds by first choosing $\kappa$ sufficiently large to conclude that there is a constant $C>0$ (independent of $k$) such that
\begin{align*}
  \E\Big[\Big(\sum_{j=1}^{N^*\wedge k} \ell(\tau^*_j,X^{j-1}_{\tau^*_{j}},\beta^*_j)\Big)^2\Big]&\leq C(1+|x|^{2\rho})
\end{align*}
and then taking the limit as $k\to\infty$, while applying the monotone convergence theorem.\\

\textbf{Step 3:} To establish \eqref{ekv:trunk-game-value-upper}, we show that for any $u\in\bar\mcU_s$ it holds that
\begin{align}\label{ekv:is-saddle-point-again}
  P^{t,x;u,\nu^{n,*}(u)}_s\leq P^{t,x;u^*,\nu^{n,*}(u^*)}_s .
\end{align}
Suppose that $\hat u\in \bar\mcU_{s}$ is another impulse control is another impulse control, then
\begin{align*}
  \bar Y^{t,x,n}_{s'}&= \bar Y^{t,x,n}_{\hat\tau_1}+\int_{s'}^{\hat\tau_1} f^n(r,X^{t,x}_r,\bar Y^{t,x,n}_r,\bar Z^{t,x,n}_r,\bar V^{t,x,n}_{r})dr - \int_{s'}^{\hat\tau_1}\bar Z^{t,x,n}_r dW_r
  \\
  &\quad-\int_{{s'}}^{\hat\tau_1}\!\!\!\int_E \bar V^{t,x,n}_{r}(e)\mu(dr,de)+\bar K^{+,t,x,n}_{\hat\tau_1}-\bar K^{+,t,x,k,n}_{s'}
  \\
  &\geq \ett_{[\hat\tau_1<T]}\{v_{k-1,n}(\hat\tau_1,X^{t,x}_{\hat\tau_1}+\xi(\hat\tau_1,X^{t,x}_{\hat\tau_1},\hat \beta_1)) -\ell({\hat\tau_1},X^{t,x}_{\hat\tau_1},\hat\beta_1)\} +\ett_{[\hat\tau_1=T]}\psi(X^{t,x}_T)
  \\
  &\quad +\int_{s'}^{\hat\tau_1} f^n(r,X^{t,x}_r,Y^{t,x,k,n}_r,Z^{t,x,k,n}_r,V^{t,x,k,n}_{r})dr - \int_{s'}^{\hat\tau_1}Z^{t,x,k,n}_r dW_r-\int_{{s'}}^{\hat\tau_1}\!\!\!\int_E V^{t,x,k,n}_{r}(e)\mu(dr,de)\\
  &\quad+\bar K^{+,t,x,k,n}_{\hat\tau_1}-\bar K^{t,x,k,n}_{s'}
\end{align*}
for all ${s'}\in [t,\hat\tau_1]$. Arguing as above gives that there is a sequence $(\hat\mcY^j,\hat\mcZ^j,\hat\mcV^j,\hat\mcK^j)_{j=0}^k\subset\mcS^2\times\mcH^2(W)\times\mcH^2(\mu)\times\mcA^2$ such that letting $\hat\mcY:=\ett_{[t,\hat\tau_{1}]}\hat\mcY^0+\sum_{j=1}^k\ett_{(\hat\tau_{j},\hat\tau_{j+1}]}\hat\mcY^j$, $\hat\mcZ:=\sum_{j=0}^k\ett_{(\hat\tau_{j},\hat\tau_{j+1}]}\hat\mcZ^j$, $\hat\mcV:=\sum_{j=0}^k\ett_{(\hat\tau_{j},\hat\tau_{j+1}]}\hat\mcV^j$ and $\hat\mcK_{s'}:=\sum_{j=0}^{k}\ett_{[\hat \tau_j<{s'}]}\{\hat\mcK^j_{{s'}\wedge\hat\tau_{j+1}}-\hat\mcK^j_{\hat\tau_j}\}$, with $\hat\tau_0:=t$, implies that $(\hat\mcY,\hat\mcZ,\hat\mcV,\hat\mcK)$ satisfies
\begin{align*}
  \hat\mcY_{s'}& \geq \supOP \bar Y^{n}(\hat\tau_{k+1},X^{t,x;[\hat u]_{k}}_{\hat\tau_{k+1}})+\int_{s'}^{\hat\tau_{k+1}} f^n\big(r,X^{t,x;\hat u}_r,\hat\mcY_r,\hat\mcZ_r,\hat\mcV_r\big)dr-\int_{s'}^{\hat\tau_{k+1}}\hat\mcZ_rdW_r
  \\
  &\quad - \int_{{s'}}^{\hat\tau_{k+1}}\!\!\!\int_E \hat\mcV_{r}(e)\mu(dr,de)
  \\
  &\quad-\sum_{j=1}^{\hat N\wedge k}\ett_{[{s'}\leq\hat\tau_j]}\ell({\hat \tau_j},X^{t,x;[\hat u]_{j-1}}_{\hat \tau_j},\hat \beta_j)+\hat\mcK_{\hat\tau_{k+1}}-\hat \mcK_{{s'}},\quad\forall {s'}\in[t,\hat\tau_{k+1}]
\end{align*}
for each $k\in\bbN$ and $\hat\mcY_{s}=\bar Y^{t,x,n}_{s}$. Now, since $\hat u\in\bar\mcU_{s}$, we have by definition that
\begin{align*}
  \E\Big[\Big(\sum_{j=1}^{\hat N}\ell({\hat\tau_j},X^{t,x;[\hat u]_{j-1}}_{\hat\tau_j},\hat\beta_j)\Big)^2\Big]<\infty
\end{align*}
and we can take the limit as $k\to\infty$ to find that
\begin{align*}
  \hat\mcY_{s'}& \geq \psi(X^{t,x;\hat u}_{T})+\int_{s'}^{T} f^n\big(r,X^{t,x;\hat u}_r,\hat\mcY_r,\hat\mcZ_r,\hat\mcV_r\big)dr-\int_{s'}^{T}\hat\mcZ_rdW_r
  \\
  &\quad - \int_{{s'}}^{T}\!\!\!\int_E \hat\mcV_{r}(e)\mu(dr,de)
  \\
  &\quad-\sum_{j=1}^{\hat N}\ett_{[{s'}\leq\hat\tau_j]}\ell({\hat \tau_j},X^{t,x;[\hat u]_{j-1}}_{\hat \tau_j},\hat \beta_j)+\hat\mcK_{T}-\hat \mcK_{{s'}},\quad\forall {s'}\in[t,T].
\end{align*}
Comparison then gives that $P^{t,x;\hat u,\nu^{n,*}(\hat u)}_{s}\leq\hat\mcY_{s}$ and we conclude that $P^{t,x;\hat u,\nu^{n,*}(\hat u)}_{s} \leq \bar Y^{t,x,n}_s$, effectively showing that \eqref{ekv:is-saddle-point-again} holds.\qed\\

By iterating a comparison result for standard BSDEs we conclude that the sequence $(\bar v_n)_{n\in\bbN}$ is non-increasing. We conclude that there is a map $\bar v:[0,T]\times \R\to \bar \R$ such that $\bar v_n\searrow \bar v$ pointwisely and since each $\bar v_n$ is continuous it follows that $\bar v$ is u.s.c.

\subsection{Viscosity properties of $\underline v$ and $\bar v$}
By Propositions \ref{prop:dblunder} and \ref{prop:trunk-game-value-upper}, we have
\begin{align}
  \underline v(s,X^{t,x}_s)\leq \esssup_{u\in\mcU^{t,x}_s}\essinf_{\nu\in\mcV^{t,x}}P_s^{t,x;u,\nu(u)} \leq \essinf_{\nu\in\mcV^{t,x}}\esssup_{u\in\mcU^{t,x}_s}P_s^{t,x;u,\nu(u)}\leq \bar v(s,X^{t,x}_s)\label{ekv:game-value-squeeze}
\end{align}
and, in particular, $\underline v\leq V^-\leq V^+\leq \bar v$. Moreover, $\underline v_0\leq \underline v\leq \bar v\leq \bar v_0$ and we conclude that $\underline v,\bar v\in\Pi^g$. Another consequence is that, $\underline v$ and $\bar v$ are both continuous at $t=T$, with
\begin{align}
  \underline v(T,x)=\bar v(T,x)=\psi(x).\label{ekv:terminal-holds}
\end{align}

\begin{prop}\label{prop:solves-pde}
The lower and the upper limit coincide and $v=\underline v=\bar v\in\Pi^g_c$ is a the unique viscosity solution to \eqref{ekv:var-ineq} in $\Pi^g$.
\end{prop}

\noindent\emph{Proof.} We demonstrate that $\underline v$ is a supersolution and $\bar v$ is a subsolution to \eqref{ekv:var-ineq}. By the comparison result for solutions of \eqref{ekv:var-ineq}, reported in Proposition~\ref{app:prop:comp-visc}, it follows that $\bar v \leq \underline v$ and \eqref{ekv:game-value-squeeze} then allows us to conclude that $\underline v=\bar v$. In particular, we find that $v=\underline v=\bar v$ is continuous and a viscosity solution to \eqref{ekv:var-ineq}.

By continuity of $\underline v_k$ and $\bar v_n$, we have (see \eg \cite{Barles94}, p. 91),
\begin{align*}
  \underline v(t,x)=\underline v_*(t,x)&=\liminf_{k\to\infty}\!{}_*\underline v_k(t,x):=\liminf_{(k,t',x')\to(\infty,t,x) }\underline v_k(t',x'),
  \\
   \bar v(t,x)=\bar v^*(t,x)&=\limsup_{n\to\infty}\!{}^*\bar v_n(t,x):=\limsup_{(n,t',x')\to(\infty,t,x) }\bar v_n(t',x').
\end{align*}

That the terminal conditions hold for both $\underline v$ and $\bar v$ is immediate from \eqref{ekv:terminal-holds}. The remaining viscosity properties are examined below:\\

\emph{Supersolution property of $\underline v$.} For $(t,x)\in [0,T)\times\R^d$ and $(p,q,M)\in\bar J^-\underline v(t,x)$ there is (by Lemma 6.1 in~\cite{UsersGuide}) a sequence $(t_j,x_j)_{j\geq 0}$ in $[0,T)\times\R^d$ and sequences $k_j\to\infty$ and $(p_j,q_j,M_j)_{j\geq 0}$ with $(p_j,q_j,M_j)\in J^-\underline v_{k_j}(t_j,x_j)$ such that
\begin{align*}
  (t_j,x_j,\underline v_{k_j}(t_j,x_j),p_j,q_j,M_j)\to (t,x,\underline v(t,x),p,q,M).
\end{align*}
As $(p_j,q_j,M_j)\in J^-\underline v_{k_j}(t_j,x_j)$ it holds that
\begin{align*}
  &\min\{\underline v_{k_j}(t_j,x_j)-\supOP \underline v_{k_j-1}(t_j,x_j),\max\{\underline v_{k_j}(t_j,x_j)-\infOP \underline v_{k_j}(t_j,x_j),-p_j-<b(t_j,x_j),q_j>
  \\&-\frac{1}{2}\trace(\sigma\sigma^\top(t_j,x_j)M_j)-f(t_j,x_j,\underline v_{k_j}(t_j,x_j),\sigma^\top(t_j,x_j)q_j)\}\}\geq 0.
\end{align*}
We first note that
\begin{align*}
\underline v(t,x)&=\lim_{k\to\infty}\underline v_k(t,x)\geq \lim_{k\to\infty}\supOP\underline v_k(t,x)=\lim_{k\to\infty}\sup_{b\in U}\{\underline v_k(t,x+\xi(t,x,b))-\ell(t,x,b)\}
\\
&\geq\sup_{b\in U}\lim_{k\to\infty}\{\underline v_k(t,x+\xi(t,x,b))-\ell(t,x,b)\}=\supOP\underline v(t,x).
\end{align*}
It thus remains to show that
\begin{align*}
  \max\{\underline v(t,x)-\infOP \underline v(t,x),-p-<b(t,x),q> -\frac{1}{2}\trace(\sigma\sigma^\top(t,x)M)-f(t,x,\underline v(t,x),\sigma^\top(t,x)q)&\}\geq 0.
\end{align*}
Consider the case when $\underline v(t,x)<\infOP \underline v(t,x)$. Since $(\infOP\underline v_k)_{k\in\bbN}$ is a non-decreasing sequence of continuous functions (see Remark~\ref{rem:mcM-monotone}) such that $\infOP \underline v_k\nearrow \infOP \underline v$, we have (by again appealing to \cite{Barles94}, p. 91) that
\begin{align*}
\infOP\underline v(t,x)=\liminf_{k\to\infty}\!{}_*\infOP\underline v_k(t,x)
\end{align*}
and we immediately get that
\begin{align*}
\liminf_{j\to\infty} \infOP \underline v_{k_j}(t_j,x_j) \geq \infOP\underline v(t,x).
\end{align*}
In particular, we find that for $j$ sufficiently large, $\underline v_{k_j}(t_j,x_j) < \infOP \underline v_{k_j}(t_j,x_j)$. Hence,
\begin{align*}
  &-p_j-<b(t_j,x_j),q_j>-\frac{1}{2}\trace(\sigma\sigma^\top(t_j,x_j)M_j)-f(t_j,x_j,\underline v_{k_j}(t_j,x_j),\sigma^\top(t_j,x_j)q_j)\geq 0
\end{align*}
for $j$ large enough and by taking the limit as $j\to\infty$ it follows that
\begin{align*}
 -p-<b(t,x),q> -\frac{1}{2}\trace(\sigma\sigma^\top(t,x)M)-f(t,x,\underline v(t,x),\sigma^\top(t,x)q)\geq 0.
\end{align*}
Since this holds whenever $\underline v(t,x)<\infOP \underline v(t,x)$ the supersolution property of $\underline v$ follows.\\

\emph{Subsolution property of $\bar v$.} For $(t,x)\in [0,T)\times\R^d$ and $(p,q,M)\in\bar J^+\bar v(t,x)$ there is (again by Lemma 6.1 in~\cite{UsersGuide}) a sequence $(t_j,x_j)_{j\geq 0}$ in $[0,T)\times\R^d$ and sequences $n_j\to\infty$ and $(p_j,q_j,M_j)_{j\geq 0}$ with $(p_j,q_j,M_j)\in J^+\bar v_{n_j}(t_j,x_j)$ such that
\begin{align*}
  (t_j,x_j,\bar v_{n_j}(t_j,x_j),p_j,q_j,M_j)\to (t,x,\bar v(t,x),p,q,M).
\end{align*}
As $(p_j,q_j,M_j)\in J^+\bar v_{n_j}(t_j,x_j)$, it holds that
\begin{align*}
  &\min\{\bar v_{n_j}(t_j,x_j)-\supOP \bar v_{n_j}(t_j,x_j),-p_j-<b(t_j,x_j),q_j>-\frac{1}{2}\trace(\sigma\sigma^\top(t_j,x_j)M_j)
  \\& + n_j\int_E(\bar v_{n_j}(t_j,x_j+\gamma(t_j,x_j,e))+\chi(t_j,x_j,e)-\bar v_{n_j}(t_j,x_j))^-\lambda(de)
  -f(t_j,x_j,\bar v_{n_j}(t_j,x_j),\sigma^\top(t_j,x_j)q_j)\}\leq 0.
\end{align*}
First note that since $\supOP\bar v_n$ is continuous, we have
\begin{align*}
  \limsup_{j\to\infty}\supOP \bar v_{n_j}(t_j,x_j) &\leq \limsup_{n\to\infty}\!{}^*\supOP \bar v_{n}(t,x) = \supOP \bar v(t,x).
\end{align*}
When $\bar v(t,x)>\supOP \bar v(t,x)$ it thus follows that $\bar v_{n_j}(t_j,x_j)>\supOP \bar v(t_j,x_j)$ whenever $j$ is sufficiently large. Hence,
\begin{align}\nonumber
  &-p_j-<b(t_j,x_j),q_j>-\frac{1}{2}\trace(\sigma\sigma^\top(t_j,x_j)M_j)
  \\& +n_j\int_E(\bar v_{n_j}(t_j,x_j+\gamma(t_j,x_j,e))+\chi(t_j,x_j,e)-\bar v_{n_j}(t_j,x_j))^-\lambda(de)
  -f(t_j,x_j,\bar v_{n_j}(t_j,x_j),\sigma^\top(t_j,x_j)q_j)\leq 0\label{ekv:subsolu}
\end{align}
for large $j$. Sending $j\to\infty$ gives that
\begin{align*}
  -p - <b(t,x),q>-\frac{1}{2}\trace(\sigma\sigma^\top(t,x)M)-f(t,x,\bar v(t,x),\sigma^\top(t,x)q)\leq 0.
\end{align*}
Now, assume that $\bar v(t,x)>\infOP \bar v(t,x)$, in which case there exists an $e_0\in E$ such that
\begin{align*}
  (\bar v(t,x+\gamma(t,x,e_0))-\chi(t,x,e_0))-\bar v(t,x)<0.
\end{align*}
Since $\bar v$ is u.s.c., this in turn implies the existence of an $\eps>0$ and an open neighborhood $E_0\in\mcB(E)$ of $e_0$ such that
\begin{align*}
  \sup_{e\in E_0}\{\bar v(t,x+\gamma(t,x,e))-\chi(t,x,e)\}-\bar v(t,x)\leq -2\eps
\end{align*}
Now, as above we have
\begin{align*}
   \limsup_{n\to\infty}\!{}^*\sup_{e\in \bar E_0}\{\bar v_n(t,x+\gamma(t,x,e))-\chi(t,x,e)\} = \sup_{e\in \bar E_0}\{\bar v(t,x+\gamma(t,x,e))-\chi(t,x,e)\}
\end{align*}
and we conclude that
\begin{align*}
  (\bar v_{n_j}(t_j,x_j+\gamma(t,x_j,e))-\chi(t,x_j,e))-\bar v_{n_j}(t_j,x_j)\leq -\eps
\end{align*}
for all $e\in E_0$ and all $j$ sufficiently large. Consequently,
\begin{align*}
  \int_E(\bar v_{n_j}(t_j,x_j+\gamma(t_j,x_j,e))+\chi(t_j,x_j,e)-\bar v_{n_j}(t_j,x_j))^-\lambda(de)\geq \eps\lambda(E_0)
\end{align*}
for $j$ sufficiently large. However, as $\lambda(E_0)>0$, this contradicts the fact that \eqref{ekv:subsolu} holds for all $j$ and we conclude that $\bar v(t,x)\leq\infOP \bar v(t,x)$.\qed\\

The proof finishes our proof of the first part of Theorem~\ref{thm:game-value}, summarized as:
\begin{cor}
The game has a value,
\begin{align}
   \esssup_{u\in\mcU^{t,x}_s}\essinf_{\nu\in\mcV^{t,x}}P_s^{t,x;u,\nu(u)} = \essinf_{\nu\in\mcV^{t,x}}\esssup_{u\in\mcU^{t,x}_s}P_s^{t,x;u,\nu(u)} = v(s,X^{t,x}_s)
\end{align}
for each $(t,x)\in [0,T]\times\R^d$ and $s\in [t,T]$, where $v$ is the unique viscosity solution to \eqref{ekv:var-ineq} in $\Pi^g$.
\end{cor}

\section{Probabilistic representation\label{sec:prob-rep}}
What remains is to prove the probabilistic representation of $v$ in terms of solutions to \eqref{ekv:Snell-env}-\eqref{ekv:fwd-sde}. We introduce the process $Y^{t,x}_s=\esssup_{\tau\in\mcT_s}Y^{t,x;\tau}_s$, where for each $\tau\in\mcT$, the quadruple $(Y^{t,x;\tau},Z^{t,x;\tau},V^{t,x;\tau},K^{-,t,x;\tau}) \in\mcS^2_{[0,\tau]}\times \mcH^2_{[0,\tau]}(W) \times \mcH^2_{[0,\tau]}(\mu) \times \mcA^2_{[0,\tau]}$ is the unique maximal solution to the BSDE
\begin{align}\label{ekv:stopped-bsde-v}
  \begin{cases}
     Y^{t,x;\tau}_s=\supOP v(\tau,X^{t,x}_\tau)+\int_s^\tau f(r,X^{t,x}_r, Y^{t,x;\tau}_r, Z^{t,x;\tau}_r)dr -\int_s^\tau  Z^{t,x;\tau}_r dW_r
    \\
    \quad-\int_{s}^\tau\!\!\!\int_E  V^{t,x;\tau}_{r}(e)\mu(dr,de)-(K^{-,t,x;\tau}_\tau- K^{-,t,x;\tau}_s),\quad \forall s\in [0,\tau],
    \\
    V^{t,x;\tau}_s(e)\geq - \chi(s,X^{t,x}_{s-},e),\quad d\Prob\otimes ds\otimes \lambda(de)-a.e.
  \end{cases}
\end{align}
We need to show that $Y^{t,x}_s=v(s,X^{t,x}_s)$ for all $s\in [t,T]$. First, since $v\geq\underline v_{k}$ for each $k\in\bbN$, comparison together with \eqref{ekv:stopped-bsde-k} gives that $Y^{t,x}_s\geq \underline v_k(s,X^{t,x}_s)$ and thus $Y^{t,x}_s\geq v(s,X^{t,x}_s)$ for all $s\in [t,T]$.

On the other hand, if we for each $n\in\bbN$, let the quadruple $(\tilde Y^{t,x,n},\tilde Z^{t,x,n},\tilde V^{t,x,n},\tilde K^{+,t,x,n})\in\mcS^{2}_t \times \mcH^{2}_t(W) \times \mcH^{2}_t(\mu) \times \mcA^{2}_t$ be the unique solution to the reflected BSDE
\begin{align}\label{ekv:rbsde-jmp-v}
  \begin{cases}
     \tilde Y^{t,x,n}_s=\psi(X^{t,x}_T)+\int_s^T f^n(r,X^{t,x}_r, \tilde Y^{t,x,n}_r, \tilde Z^{t,x,n}_r, \tilde V^{t,x,n}_{r})dr -\int_s^T  \tilde Z^{t,x,n}_r dW_r-\int_{s}^T\!\!\!\int_E \tilde V^{t,x,n}_{r}(e)\mu(dr,de)
    \\
    \quad+\tilde K^{+,t,x,n}_T - \tilde K^{+,t,x,n}_s,\quad \forall s\in [t,T],
    \\
    \tilde Y^{t,x,n}_s\geq \supOP v(s,X^{t,x}_s),\, \forall s\in [t,T] \quad\text{and}\quad\int_t^T\!\! \big(\tilde Y^{t,x,n}_s-\supOP v(s,X^{t,x}_s)\big)dK^{+,t,x,n}_s,
  \end{cases}
\end{align}
then comparison for BSDEs with jumps along with the fact that $v\leq \bar v_n$ for each $n\in\bbN$ gives that $\tilde Y^{t,x,n}\leq \bar Y^{t,x,n}$. By Proposition~\ref{prop:rbsde-jmp} we also have the representation
\begin{align*}
  \tilde Y^{t,x,n}_s=\esssup_{\tau\in\mcT_s}\tilde Y^{t,x,n;\tau}_s,
\end{align*}
where for each $\tau\in\mcT$, the triple $(\tilde Y^{t,x,n;\tau},\tilde Z^{t,x,n;\tau},\tilde V^{t,x,n;\tau}) \in\mcS^{2}_{[0,\tau]} \times \mcH^{2}_{[0,\tau]}(W) \times \mcH^{2}_{[0,\tau]}(\mu)$ satisfies
\begin{align*}
  \tilde Y^{t,x,n;\tau}_s&=\supOP v(\tau,X^{t,x}_\tau)+\int_s^\tau f^n(r,X^{t,x}_r,\tilde Y^{t,x,n;\tau}_r,\tilde Z^{t,x,n;\tau}_r,\tilde V^{t,x,n;\tau}_{r})dr
  \\
  &\quad -\int_s^\tau \tilde Z^{t,x,n;\tau}_r dW_r-\int_{s}^\tau\!\!\!\int_E \tilde V^{t,x,n;\tau}_{r}(e)\mu(dr,de),\quad \forall s\in [0,\tau].
\end{align*}
Now, a straightforward generalization of Lemma 3.2 in \cite{Kharroubi2010} to a random horizon setting gives that $Y^{t,x;\tau}_s\leq\tilde Y^{t,x,n;\tau}_s$ for each $\tau\in\mcT_t$ and $s\in[t,\tau]$, from which we conclude that $Y^{t,x}_s\leq\tilde Y^{t,x,n}_s\leq \bar Y^{t,x,n}_s=\bar v_n(s,X^{t,x}_s)$. Since $\bar v_n\searrow v$, we thus get that $Y^{t,x}_s\leq v(s,X^{t,x}_s)$, proving that $Y^{t,x}_s= v(s,X^{t,x}_s)$.

Uniqueness in $\bigS^2$ is immediate from uniqueness of solutions to \eqref{ekv:var-ineq} and Theorem~\ref{thm:qvi-stop}.

%\begin{thm}
%$v(s,X^{t,x}_s)=Y^{t,x}_s$, where $Y^{t,x}\in\mcS^2$ satisfies $Y^{t,x}_s=\esssup_{\tau\in\mcT_s}Y^{t,x,\tau}_s$, where\\ $(Y^{t,x,\tau},Z^{t,x,\tau},V^{t,x,\tau},K^{-,t,x,\tau})\in\mcS^2_{\tau}\times\mcH^2_\tau(W)\times\mcH^2_\tau(\mu)\times\mcA^2_{\tau}$ is the unique maximal solution to
%\begin{align}\label{ekv:rbsde-final}
%  \begin{cases}
%     Y^{t,x,\tau}_s=\supOP Y(\tau,X^{t,x}_\tau)+\int_s^\tau f(r,X^{t,x}_r, Y^{t,x,\tau}_r, Z^{t,x,\tau}_r)dr-\int_s^\tau  Z^{t,x,\tau}_r dW_r-\int_{s}^\tau\!\!\!\int_E  V^{t,x,\tau}_{r}(e)\mu(dr,de)
%    \\
%    \quad-(K^{-,t,x,\tau}_\tau- K^{-,t,x,\tau}_s),\quad \forall s\in [t,\tau],
%    \\
%    V^{t,x,\tau}_s(e)\geq - \chi(s,X^{t,x}_{s-},e),\quad d\Prob\otimes ds\otimes \lambda(de)-a.e.
%  \end{cases}
%\end{align}
%\end{thm}

\appendix
\section{Comparison principle\label{app:uniq}}
We prove a comparison result for \eqref{ekv:var-ineq}, stating that a supersolution to \eqref{ekv:var-ineq} always dominates a subsolution. The proof of this principle is based on a sequence of lemmas.

\begin{lem}\label{lem:is-super}
Let $v$ be a supersolution to \eqref{ekv:var-ineq} satisfying
\begin{align*}
\forall (t,x)\in[0,T]\times\R^d,\quad |v(t,x)|\leq C(1+|x|^{2\varrho})
\end{align*}
for some $C,\varrho>0$. Then there is a $\varpi_0 > 0$ such that for any $\varpi>\varpi_0$ and $\theta > 0$, the function $(t,x)\mapsto w(t,x):=v(t,x) + \theta e^{-\varpi t}(1+((|x|-K_{\gamma,\xi})^+)^{2\varrho+2})$ is also a supersolution to \eqref{ekv:var-ineq}.
\end{lem}
%\todo[inline]{$\kappa$ used at several places}

\noindent \emph{Proof.} Since $v$ is a supersolution and $\theta e^{-\varpi T}(1+((|x|-K_{\gamma,\xi})^+)^{2\varrho+2})\geq 0$, we have $w(T,x)\geq v(T,x)\geq \psi(x)$ so that the terminal condition holds. Moreover,
\begin{align*}
  v(t,x) - \sup_{b\in U}\{v(t,x+\xi(t,x,b))-\ell(t,x,b)\}\geq 0
\end{align*}
and we have,
\begin{align*}
&w(t,x) - \sup_{b\in U}\{w(t,x+\xi(t,x,b))-\ell(t,x,b)\}
\\
&=v(t,x) + \theta e^{-\varpi t}(1+((|x|-K_{\gamma,\xi})^+)^{2\varrho+2})
\\
&\quad - \sup_{b\in U}\{v(t,x+\xi(t,x,b)) + \theta e^{-\varpi t}(1+((|x+\xi(t,x,b)|-K_{\gamma,\xi})^+)^{2\varrho+2}) - \ell(t,x,b)\}
\\
&\geq v(t,x) - \sup_{b\in U}\{v(t,x+\xi(t,x,b))-\ell(t,x,b)\}
\\
&\quad+\theta e^{-\varpi t}\{1+((|x|-K_{\gamma,\xi})^+)^{2\varrho+2}- \sup_{b\in U}(1+((|x+\xi(t,x,b)|-K_{\gamma,\xi})^+)^{2\varrho+2})\}.
\end{align*}
Now, either $|x|\leq K_{\gamma,\xi}$ in which case it follows by \eqref{ekv:imp-bound} that $|x+\xi(t,x,b)|\leq K_{\gamma,\xi}$ or $|x|> K_{\gamma,\xi}$ and \eqref{ekv:imp-bound} gives that $|x+\xi(t,x,b)|\leq |x|$. We conclude that
\begin{align*}
  w(t,x) - \sup_{b\in U}\{w(t,x+\xi(t,x,b))-\ell(t,x,b)\}\geq 0.
\end{align*}
Assume now that
\begin{align*}
  v(t,x) - \inf_{e\in E}\{v(t,x+\gamma(t,x,e))+\chi(t,x,e)\}\geq 0.
\end{align*}
In this case,
\begin{align*}
&w(t,x) - \inf_{e\in E}\{w(t,x+\gamma(t,x,e))+\chi(t,x)\}
\\
&=v(t,x) + \theta e^{-\varpi t}(1+((|x|-K_{\gamma,\xi})^+)^{2\varrho+2})
\\
&\quad - \inf_{e\in E}\{v(t,x+\gamma(t,x,e)) + \theta e^{-\varpi t}(1+((|x+\gamma(t,x,e)|-K_{\gamma,\xi})^+)^{2\varrho+2}) +\chi(t,x,e)\}
\\
&\geq v(t,x) - \inf_{e\in E}\{v(t,x+\gamma(t,x,e))+\chi(t,x,e)\}
\\
&\quad+\theta e^{-\varpi t}\{1+((|x|-K_{\gamma,\xi})^+)^{2\varrho+2}- \sup_{e\in E} (1+((|x+\gamma(t,x,e)|-K_{\gamma,\xi})^+)^{2\varrho+2})\}.
\end{align*}
and we can repeat the above argument to conclude that
\begin{align*}
  w(t,x) - \inf_{e\in E}\{w(t,x+\gamma(t,x,e))+\chi(t,x,e)\}\geq 0.
\end{align*}
Consider instead the case when
\begin{align*}
  v(t,x) - \inf_{e\in E}\{v(t,x+\gamma(t,x,e))+\chi(t,x,e)\}< 0
\end{align*}
and let $\varphi\in C^{1,2}([0,T]\times\R^d\to\R)$ be such that $\varphi-w$ has a local maximum of 0 at $(t,x)$ with $t<T$. Then $(\tilde t,\tilde x)\mapsto\tilde \varphi(\tilde t,\tilde x):=\varphi (\tilde t,\tilde x)-\theta e^{-\varpi \tilde t}(1+((|\tilde x|-K_{\gamma,\xi})^+)^{2\varrho+2})\in C^{1,2}([0,T]\times\R^d\to\R)$ and $\tilde \varphi-v$ has a local maximum of 0 at $(t,x)$. Since $v$ is a viscosity supersolution, we have
\begin{align*}
    &-\partial_t(\varphi(t,x)-\theta e^{-\varpi t}(1+((|x|-K_{\gamma,\xi})^+)^{2\varrho+2}))-\mcL (\varphi(t,x)-\theta e^{-\varpi t}(1+((|x|-K_{\gamma,\xi})^+)^{2\varrho+2}))
    \\
    &-\tilde f(t,x,\varphi(t,x)-\theta e^{-\varpi t}(1+((|x|-K_{\gamma,\xi})^+)^{2\varrho+2}),\sigma^\top(t,x)\nabla_x (\varphi(t,x)-\theta e^{-\varpi t}(1+((|x|-K_{\gamma,\xi})^+)^{2\varrho+2})))\geq 0.
\end{align*}
Consequently,
\begin{align*}
&-\partial_t\varphi(t,x)-\mcL \varphi(t,x)-\tilde f(t,x,\varphi(t,x),\sigma^\top(t,x)\nabla_x \varphi(t,x))
\\
&\geq \theta \varpi e^{-\varpi t}(1+((|x|-K_{\gamma,\xi})^+)^{2\varrho+2})-\theta\mcL e^{-\varpi t}(1+((|x|-K_{\gamma,\xi})^+)^{2\varrho+2})
\\
&\quad \tilde f(t,x,\varphi(t,x)-\theta e^{-\varpi t}(1+((|x|-K_{\gamma,\xi})^+)^{2\varrho+2}),\sigma^\top(t,x)\nabla_x (\varphi(t,x)-\theta e^{-\varpi t}(1+((|x|-K_{\gamma,\xi})^+)^{2\varrho+2})))
\\
&\quad-\tilde f(t,x,\varphi(t,x),\sigma^\top(t,x)\nabla_x \varphi(t,x))
\\
&\geq \theta \varpi e^{-\varpi t}(1+((|x|-K_{\gamma,\xi})^+)^{2\varrho+2})-\theta C(1+\varrho) e^{-\varpi t}(1+((|x|-K_{\gamma,\xi})^+)^{2\varrho+2})
\\
&\quad -k_f(1+\varrho)\theta e^{-\varpi t}(1+((|x|-K_{\gamma,\xi})^+)^{2\varrho+2}),
\end{align*}
where the right hand side is non-negative for all $\theta> 0$ and all $\varpi>\varpi_0$ for some $\varpi_0>0$.\qed\\

We have the following result, the proof of which we omit since it is a straightforward generalization of a classical result:
\begin{lem}\label{lem:integ-factor}
For any $\kappa\in\R$, a locally bounded function $v:[0,T]\times \R^d\to\R$ is a viscosity supersolution (resp. subsolution) to \eqref{ekv:var-ineq} if and only if $\check v(t,x):=e^{\kappa t}v(t,x)$ is a viscosity supersolution (resp. subsolution) to
\begin{align}\label{ekv:var-ineq-if}
\begin{cases}
  \min\{\check v(t,x)-\sup_{b\in U}\{\check v(t,x+\xi(t,x,b))-e^{\kappa t}\ell(t,x,b)\},\max\{\check v(t,x)-
  \\
\inf_{e\in E}\{\check v(t,x+\gamma(t,x,e))+e^{\kappa t}\chi(t,x,e)\},-\check v_t(t,x)+\kappa \check v(t,x)
  \\
-\mcL \check v(t,x)-e^{\kappa t}\tilde f(t,x,e^{-\kappa t}\check v(t,x),e^{-\kappa t}\sigma^\top(t,x)\nabla_x \check v(t,x))\}=0,\quad\forall (t,x)\in[0,T)\times \R^d \\
  \check v(T,x)=e^{\kappa T}\psi(x).
\end{cases}
\end{align}
\end{lem}
\begin{rem}
Here, it is important to note that $\check \ell(t,x,b):=e^{\kappa t}\ell(t,x,b)$, $\check \chi(t,x,e):=e^{\kappa t}\chi(t,x,e)$, $\check f(t,x,y,z):=-\kappa y+e^{\kappa t} f(t,x,e^{-\kappa t}y,e^{-\kappa t}z)$ and $\check \psi(x):=e^{\kappa T}\psi(x)$ satisfy Assumption~\ref{ass:oncoeff}. In particular, this implies that Lemma~\ref{lem:is-super} holds for supersolutions to \eqref{ekv:var-ineq-if} as well. We thus define the maps $\check\supOP v(t,x):=\sup_{b\in U}\{v(t,x+\xi(t,x,b))-\check\ell(t,x,b)\}$ and $\check\infOP v(t,x):=\inf_{e\in E}\{v(t,x+\gamma(t,x,e))+\check\chi(t,x,e)\}$.
\end{rem}

The above lemma allows us to consider super- and sub-solutions to \eqref{ekv:var-ineq-if} rather than those of \eqref{ekv:var-ineq}. Assume that $v$ (resp.~$u$) is a supersolution (resp. subsolution) to \eqref{ekv:var-ineq-if}. Since $v$ is l.s.c.~and $u$ is u.s.c., $u - v$ is u.s.c.~implying that there is, for each $R > 0$, a point $(t_0,x_0)\in [0,T]\times\bar B_R$ (where $\bar B_R$ is the closed ball of radius $R$ in $\R^d$, centered at the origin) such that
\begin{align*}
  \max_{(t,x)\in [0,T] \times \bar B_R} \{u(t,x)-v(t,x)\}= u(t_0,x_0)-v(t_0,x_0).
\end{align*}
We let $A_{\max,R}(v,u):=\{(t,x)\in [0,T]\times \bar B_R:u(t,x)-v(t,x)=\max_{(t',x')\in [0,T]\times\bar B_R} \{u(t',x')-v(t',x')\}\}$.

\begin{lem}\label{app:lem:off-boundaries}
Assume that $v$ is a supersolution to \eqref{ekv:var-ineq-if} while $u$ is a subsolution. If $\max_{(t',x')\in [0,T]\times\bar B_R} \{u(t',x')-v(t',x')\}>0$ for some $R\geq K_{\gamma,\xi}$, then there is a point $(t,x)\in A_{\max,R}(v,u)$ such that $v(t,x) < \check\infOP v(t,x)$ and $u(t,x) > \check\supOP u(t,x)$.
\end{lem}

\noindent\emph{Proof.} We let $\eta:=\max_{(t,x)\in [0,T]\times\bar B_R} \{u(t,x)-v(t,x)\}$. We search for a contradiction and assume that either $v(t,x) \geq \check\infOP v(t,x)$ or $u(t,x) \leq \check\supOP u(t,x)$ whenever $(t,x)\in A_{\max,R}(v,u)$. We then pick an arbitrary $(t,x_0)\in A_{\max,R}(v,u)$ and note that $t\in [0,T)$. If $u(t,x_0) \leq \check\supOP u(t,x_0)$, then by upper semi-continuity, there is a $b_1\in U$ such that
\begin{align}\label{ekv:since-sub}
  u(t,x_0) \leq u(t,x_0+\xi(t,x_{0},b_1))-\check\ell(t,x_0,b_1) = u(t,x_1)-\check\ell(t,x_0,b_1),
\end{align}
with $x_1 := x_0+\xi(t,x_{0},b_1)\in\bar B_R$. Hence,
\begin{align*}
  \eta = u(t,x_0) - v(t,x_0) \leq u(t,x_1)-\check\ell(t,x_0,b_1) - v(t,x_0).
\end{align*}
On the other hand, since $v$ is a supersolution,
\begin{align}\label{ekv:since-super}
v(t,x_0)\geq v(t,x_1)-\check\ell(t,x_0,b_1)
\end{align}
and we conclude that
\begin{align*}
  \eta \leq u(t,x_1)-\check\ell(t,x_0,b_1) - v(t,x_0)\leq u(t,x_1)-\check\ell(t,x_0,b_1) - (v(t,x_1)-\check\ell(t,x_0,b_1)).
\end{align*}
In particular, $(t,x_1)\in A_{\max,R}(v,u)$ and equality must hold in both \eqref{ekv:since-sub} and \eqref{ekv:since-super}. If instead $u(t,x_0) \geq \check\supOP u(t,x)$, our assumption implies that $v(t,x_1) \geq v(t,x_1)+\check\chi(t,x_0,e_1)$ for some $e_1\in E$, with $x_1 := x_0+\gamma(t,x_{0},e_1)\in\bar B_R$. Then we can repeat the above argument to conclude that $u(t,x_1) = u(t,x_1)+\check\chi(t,x_0,e_1)$, $v(t,x_1) = v(t,x_1)+\check\chi(t,x_0,e_1)$ and $(t,x_1)\in A_{\max,R}(v,u)$. Now, as $(t,x_1)\in A_{\max,R}(v,u)$ we can repeat the process to find a $b_2\in U$ (or an $e_2\in E$) such that,  for $\phi=u,v$, it holds that $\phi(t,x_1) = \phi(t,x_2)+\check\ell(t,x_0,b_2)$, with $x_2 := x_1+\xi(t,x_{1},b_2)$ (or $\phi(t,x_1) = \phi(t,x_2) - \check\chi(t,x_0,e_2)$, with $x_2 := x_1+\gamma(t,x_{1},e_2)\in\bar B_R$). This process can be repeated indefinitely to find a sequence of points $(x_i,b_i,e_i,\iota_i)_{i\in\bbN}$ in $\bar B_R\times U\times E\times \{1,2\}$ such that for each $i\in\bbN$,
\begin{align*}
  x_{i+1}=x_{i}+\ett_{[\iota_{i+1}=1]}\xi(t,x_{i},b_{i+1})+\ett_{[\iota_{i+1}=2]}\gamma(t,x_{i},e_{i+1})
\end{align*}
and
\begin{align*}
  \phi(t,x_i)=\phi(t,x_{i+1})-\ett_{[\iota_{i+1}=1]}\check\ell(t,x_i,b_{i+1})+\ett_{[\iota_{i+1}=2]}\check\chi(t,x_i,e_{i+1})
\end{align*}
for $\phi=u,v$. Since $\bar B_{R}$ is compact, there is a $\hat x\in \bar B_{R}$ and an increasing sequence of integers $(\pi_i)_{i\in\bbN}$ such that $x_{\pi_i}\to \hat x$. Moreover, polynomial growth implies that $v$ is bounded on $\bar B_{R}$ and it follows that there is a subsequence $(\pi'_i)_{i\in\bbN}$ such that $\lim_{i\to\infty}v(t,x_{\pi'_i})$ exists. However, for $i$ sufficiently large, $|x_{\pi'_{i+1}}-x_{\pi'_i}|\leq h_1(t)$ implying by Assumption~\ref{ass:nofreeloop} that
\begin{align*}
  |v(t,x_{\pi'_i})-v(t,x_{\pi'_{i+1}})| = \big|\sum_{j=\pi'_i+1}^{\pi'_{i+1}}(\ett_{[\iota_j=1]}\check\ell(t,x_{j-1},b_j)-\ett_{[\iota_j=2]}\check\chi(t,x_{j-1},e_j))\big|\geq h_2(t)
\end{align*}
contradicting the existence of the limit $\lim_{i\to\infty}v(t,x_{\pi'_i})$.\qed\\

Using the results presented in the above lemmas we derive the following comparison principle for viscosity solutions in $\Pi^g$:

\begin{prop}\label{app:prop:comp-visc}
Let $v$ (resp. $u$) be a supersolution (resp. subsolution) to \eqref{ekv:var-ineq}. If $u,v\in \Pi^g$, then $u\leq v$.
\end{prop}

\noindent \emph{Proof.} First, we note that it is sufficient to show that the statement holds for solutions to~\eqref{ekv:var-ineq-if} for some $\kappa\in \R$. We thus assume that $v$ (resp.~$u$) is a viscosity supersolution (resp.~subsolution) to \eqref{ekv:var-ineq-if} for $\kappa\in \R$ specified below.% Furthermore, we may without loss of generality assume that $v$ is l.s.c.~and $u$ is u.s.c.

By assumption, $u,v\in \Pi^g$, which implies that there are constants $C>0$ and $\varrho>0$ such that
\begin{align}\label{ekv:uv-bound}
|v(t,x)|+|u(t,x)|\leq C(1+|x|^{2\varrho}),\quad \forall (t,x)\in[0,T]\times\R^d.
\end{align}
Now, for any $\varpi>0$ we only need to show that
\begin{align*}
w(t,x)&=w^{\theta,\varpi}(t,x):=v(t,x)+\theta e^{-\varpi t}(1+((|x|-K_{\gamma,\xi})^+)^{2\varrho+2})
\\
&\geq u(t,x)
\end{align*}
for all $(t,x)\in[0,T]\times\R^d$ and any $\theta>0$. Then the result follows by taking the limit $\theta\to 0$. We know from Lemma~\ref{lem:is-super} that there is a $\varpi_0>0$ such that $w$ is a supersolution to \eqref{ekv:var-ineq-if} for each $\varpi\geq\varpi_0$ and $\theta>0$. We thus assume that $\varpi\geq\varpi_0$.

We search for a contradiction and assume that there is a $(t_0,x_0)\in [0,T]\times \R^d$ such that $u(t_0,x_0)>w(t_0,x_0)$. By \eqref{ekv:uv-bound}, there is for each $\theta>0$ a $R> K_{\gamma,\xi}$ such that
\begin{align*}
w(t,x)\geq u(t,x),\quad\forall (t,x)\in[0,T]\times\R^d,\:|x|\geq R.
\end{align*}
Our assumption thus implies that there is a point $(\bar t,\bar x)\in[0,T)\times B_R$ such that
\begin{align*}
\max_{(t,x)\in[0,T]\times\R^d}(u(t,x)-w(t,x))&=\max_{(t,x)\in[0,T)\times B_R}(u(t,x)-w(t,x))
\\
&=u(\bar t,\bar x)-w(\bar t,\bar x)=\eta>0.
\end{align*}
By Lemma~\ref{app:lem:off-boundaries} there is then at least one point $(t^*,x^*)\in[0,T)\times B_R$ such that
\begin{enumerate}[a)]
  \item $u(t^*,x^*)-w(t^*,x^*)= \eta$,
  \item $u(t^*,x^*)>\sup_{b\in U}\{v(t^*,x^*+\xi(t^*,x^*,b))-\check \ell(t^*,x^*,b)\}$ and
  \item $w(t^*,x^*)<\inf_{e\in E}\{w(t^*,x^*+\gamma(t^*,x^*,e))+\check \chi(t^*,x^*,e)\}$.
\end{enumerate}

The remainder of the proof follows along the lines of the proof of Proposition 4.1 in~\cite{Morlais13} (see also Proposition 6.4 in \cite{PerningeJMAA23}) and is included only for the sake of completeness.\\

We assume without loss of generality that $\varrho\geq 2$ and define
\begin{align*}
\Phi_n(t,x,y):=u(t,x)-w(t,y)-\varphi_n(t,x,y),
\end{align*}
where
\begin{align*}
  \varphi_n(t,x,y):=\frac{n}{2}|x-y|^{2\varrho}+|x-x^*|^2+|y-x^*|^2+(t-t^*)^2.
\end{align*}
Since $u$ is u.s.c.~and $w$ is l.s.c.~there is a triple $(t_n,x_n,y_n)\in[0,T]\times \bar B_R\times \bar B_R$ (with $\bar B_R$ the closure of $B_R$) such that
\begin{align*}
  \Phi_n(t_n,x_n,y_n)=\max_{(t,x,y)\in [0,T]\times \bar B_R\times \bar B_R}\Phi_n(t,x,y).
\end{align*}
Now, the inequality $2\Phi_n(t_n,x_n,y_n)\geq \Phi_n(t_n,x_n,x_n)+\Phi_n(t_n,y_n,y_n)$ gives
\begin{align*}
n|x_n-y_n|^{2\varrho}\leq u(t_n,x_n)-u(t_n,y_n)+w(t_n,x_n)-w(t_n,y_n).
\end{align*}
Consequently, $n|x_n-y_n|^{2\varrho}$ is bounded (since $u$ and $w$ are bounded on $[0,T]\times\bar B_R\times\bar B_R$) and $|x_n-y_n|\to 0$ as $n\to\infty$. We can, thus, extract subsequences $n_l$ such that $(t_{n_l},x_{n_l},y_{n_l})\to (\tilde t,\tilde x,\tilde x)$ as $l\to\infty$. Since
\begin{align*}
u(t^*,x^*)-w(t^*,x^*)\leq \Phi_n(t_n,x_n,y_n)\leq u(t_n,x_n)-w(t_n,y_n),
\end{align*}
it follows that
\begin{align*}
u(t^*,x^*)-w(t^*,x^*)&\leq \limsup_{l\to\infty} \{u(t_{n_l},x_{n_l})-w(t_{n_l},y_{n_l})\}
\\
&\leq u(\tilde t,\tilde x)-w(\tilde t,\tilde x)
\end{align*}
and as the right-hand side is dominated by $u(t^*,x^*)-w(t^*,x^*)$ we conclude that
\begin{align*}
  u(\tilde t,\tilde x)-w(\tilde t,\tilde x)=u(t^*,x^*)-w(t^*,x^*).
\end{align*}
In particular, this gives that $\lim_{l\to\infty}\Phi_n(t_{n_l},x_{n_l},y_{n_l})=u(\tilde t,\tilde x)-w(\tilde t,\tilde x)$ which implies that
\begin{align*}
  \limsup_{l\to\infty} n_l|x_{n_l}-y_{n_l}|^{2\varrho}= 0
\end{align*}
and
\begin{align*}
  (t_{n_l},x_{n_l},y_{n_l})\to (t^*,x^*,x^*).
\end{align*}
Since $u$ is u.s.c. and $w$ is l.s.c., we must also have
\begin{align*}
  (u(t_{n_l},x_{n_l}),w(t_{n_l},y_{n_l}))\to (u(t^*,x^*),w(t^*,x^*)).
\end{align*}
We can thus extract a subsequence $(\tilde n_l)_{l\geq 0}$ of $(n_l)_{l\geq 0}$ such that $t_{\tilde n_l}<T$, $|x_{\tilde n_l}|\vee |y_{\tilde n_l}|<R$ and
\begin{align*}
  u(t_{\tilde n_l},x_{\tilde n_l})-w(t_{\tilde n_l},y_{\tilde n_l})\geq \frac{\eta}{2}
\end{align*}
for all $l\in\bbN$. Moreover, since $(t,x)\mapsto \inf_{e\in E}\{w(t,x+\gamma(t,x,e))+\check \chi(t,x,e)\}$ is l.s.c.~and $(t,x)\mapsto \sup_{b\in U}\{u(t,x+\xi(t,x,b))-\check \ell(t,x,b)\}$ is u.s.c.~(see Remark~\ref{rem:mcM-monotone}), there is an $l_0\geq 0$ such that
\begin{align*}
  w(t_{\tilde n_l},y_{\tilde n_l})-\inf_{e\in E}\{w(t_{\tilde n_l},y_{\tilde n_l}+\gamma(t_{\tilde n_l},y_{\tilde n_l},e))-\check \chi(t_{\tilde n_l},y_{\tilde n_l},e)\}<0
\end{align*}
and
\begin{align*}
  u(t_{\tilde n_l},x_{\tilde n_l})-\sup_{b\in U}\{u(t_{\tilde n_l},x_{\tilde n_l}+\xi(t_{\tilde n_l},x_{\tilde n_l},b))-\check \ell(t_{\tilde n_l},x_{\tilde n_l},b)\}>0,
\end{align*}
for all $l\geq l_0$. To simplify notation we will, from now on, denote $(\tilde n_l)_{l\geq l_0}$ simply by $n$.\\

By Theorem 8.3 of~\cite{UsersGuide} there are $(p^u_n,q^u_n,M^u_n)\in \bar J^{+}u(t_n,x_n)$ and $(p^w_n,q^w_n,M^w_n)\in \bar J^{-}w(t_n,y_n)$, such that
\begin{align*}
\begin{cases}
  p^u_n-p^w_n=\partial_t\varphi_n(t_n,x_n,y_n)=2(t_n-t^*)
  \\
  q^u_n=D_x\varphi_n(t_n,x_n,y_n)=n\varrho(x-y)|x-y|^{2\varrho-2}+2(x-x^*)
  \\
  q^w_n=-D_y\varphi_n(t_n,x_n,y_n)=n\varrho(x-y)|x-y|^{2\varrho-2}+2(x-x^*)
\end{cases}
\end{align*}
and for every $\epsilon>0$,
\begin{align*}
  \left[\begin{array}{cc} M^n_x & 0 \\ 0 & -M^n_y\end{array}\right]\leq B(t_n,x_n,y_n)+\epsilon B^2(t_n,x_n,y_n),
\end{align*}
where $B(t_n,x_n,y_n):=D^2_{(x,y)}\varphi_n(t_n,x_n,y_n)$. Now, we have
\begin{align*}
  D^2_{(x,y)}\varphi_n(t,x,y)=\left[\begin{array}{cc} D_x^2\varphi_n(t,x,y) & D^2_{yx}\varphi_n(t,x,y) \\ D^2_{xy}\varphi_n(t,x,y) & D_y^2\varphi_n(t,x,y)\end{array}\right]
  = \left[\begin{array}{cc} n\xi(x,y)+2 I & -n\xi(x,y) \\ -n\xi(x,y) & n\xi(x,y) +2 I \end{array}\right]
\end{align*}
where $I$ is the identity-matrix of suitable dimension and
\begin{align*}
  \xi(x,y):=\varrho|x-y|^{2\varrho-4}\{|x-y|^2I+2(\varrho-1)(x-y)(x-y)^\top\}.
\end{align*}
In particular, since $x_n$ and $y_n$ are bounded, choosing $\epsilon:=\frac{1}{n}$ gives that
\begin{align}\label{ekv:mat-bound}
  \tilde B_n:=B(t_n,x_n,y_n)+\epsilon B^2(t_n,x_n,y_n)\leq Cn|x_n-y_n|^{2\varrho-2}\left[\begin{array}{cc} I & -I \\ -I & I \end{array}\right]+C I.
\end{align}
By the definition of viscosity supersolutions and subsolutions we have that
\begin{align*}
&-p^u_n+\kappa u(t_n,x_n)-a^\top(t_n,x_n)q^u_n-\frac{1}{2}\trace [\sigma^\top(t_n,x_n)M^u_n\sigma(t_n,x_n)]
\\
&-e^{\kappa t_n} f(t_n,x_n,e^{-\kappa t_n}u(t_n,x_n),e^{-\kappa t_n}\sigma^\top(t_n,x_n)q^u_n)\leq 0
\end{align*}
and
\begin{align*}
&-p^w_n+\kappa w(t_n,y_n)-a^\top(t_n,y_n)q^w_n-\frac{1}{2}\trace [\sigma^\top(t_n,y_n)M^w_n\sigma(t_n,y_n)]
\\
&-e^{\kappa t_n} f(t_n,y_n,e^{-\kappa t_n}w(t_n,y_n),e^{-\kappa t_n}\sigma^\top(t_n,x_n)q^w_n)\geq 0.
\end{align*}
Combined, this gives that
\begin{align*}
\kappa (u(t_n,x_n)-w(t_n,y_n))&\leq p^u_n+a^\top(t_n,x_n)q^u_n+\frac{1}{2}\trace [\sigma^\top(t_n,x_n)M^u_n\sigma(t_n,x_n)]
\\
&+e^{\kappa t_n} f(t_n,x_n,e^{-\kappa t_n}u(t_n,x_n),e^{-\kappa t_n}\sigma^\top(t_n,x_n)q^u_n)
\\
&-p^w_n-a^\top(t_n,y_n)q^w_n-\frac{1}{2}\trace [\sigma^\top(t_n,y_n)M^w_n\sigma(t_n,y_n)]
\\
&-e^{\kappa t_n} f(t_n,y_n,e^{-\kappa t_n}w(t_n,y_n),e^{-\kappa t_n}\sigma^\top(t_n,x_n)q^w_n)
\end{align*}
Collecting terms we have that
\begin{align*}
p^u_n-p^w_n&=2(t_n-t^*)
\end{align*}
and since $a$ is Lipschitz continuous in $x$ and bounded on $\bar B_R$, we have
\begin{align*}
  a^\top(t_n,x_n)q^u_n-a^\top(t_n,y_n)q^w_n&\leq  (a^\top(t_n,x_n)-a^\top(t_n,y_n))n\varrho(x_n-y_n)|x_n-y_n|^{2\varrho-2}
  \\
  &\quad+C(|x_n-x^*|+|y_n-x^*|)
  \\
  &\leq C(n|x_n-y_n|^{2\varrho}+|x_n-x^*|+|y_n-x^*|),
\end{align*}
where the right-hand side tends to 0 as $n\to\infty$. Let $s_x$ denote the $i^{\rm th}$ column of $\sigma(t_n,x_n)$ and let $s_y$ denote the $i^{\rm th}$ column of $\sigma(t_n,y_n)$ then by the Lipschitz continuity of $\sigma$ and \eqref{ekv:mat-bound}, we have
\begin{align*}
  s_x^\top M^u_n s_x-s_y^\top M^w_n s_y&=\left[\begin{array}{cc} s_x^\top  & s_y^\top \end{array}\right]\left[\begin{array}{cc} M^u_n  & 0 \\ 0 &-M^w_n\end{array}\right]\left[\begin{array}{c} s_x \\ s_y \end{array}\right]
  \\
  &\leq \left[\begin{array}{cc} s_x^\top  & s_y^\top \end{array}\right]\tilde B_n\left[\begin{array}{c} s_x \\ s_y \end{array}\right]
  \\
  &\leq C(n|x_n-y_n|^{2\varrho}+|x_n-y_n|)
\end{align*}
and we conclude that
\begin{align*}
  \limsup_{n\to\infty}\frac{1}{2}\trace [\sigma^\top(t_n,x_n)M^u_n\sigma(t_n,x_n)-\sigma^\top(t_n,y_n)M^w_n\sigma(t_n,y_n)]\leq 0.
\end{align*}
Finally, we have that
\begin{align*}
  &e^{\kappa t_n} f(t_n,x_n,e^{-\kappa t_n}u(t_n,x_n),e^{-\kappa t_n}\sigma^\top(t_n,x_n)q^u_n)-e^{\kappa t_n} f(t_n,y_n,e^{-\kappa t_n}w(t_n,y_n),e^{-\kappa t_n}\sigma^\top(t_n,x_n)q^w_n)
  \\
  &\leq k_f(u(t_n,x_n)-w(t_n,y_n)+|\sigma^\top(t_n,x_n)q^u_n-\sigma^\top(t_n,x_n)q^w_n|)
  \\
  &\quad + e^{\kappa t_n}|f(t_n,x_n,e^{-\kappa t_n}u(t_n,x_n),e^{-\kappa t_n}\sigma^\top(t_n,x_n)q^u_n)-f(t_n,y_n,e^{-\kappa t_n}u(t_n,x_n),e^{-\kappa t_n}\sigma^\top(t_n,x_n)q^u_n)|
\end{align*}
Repeating the above argument and using that $f$ is jointly continuous in $(t,x)$ uniformly in $(y,z)$ we get that the upper limit of the right-hand side when $n\to\infty$ is bounded by $k_f(u(t_n,x_n)-w(t_n,y_n))$. Put together, this gives that
\begin{align*}
(\kappa-k_f) \limsup_{n\to\infty}(u(t_n,x_n)-w(t_n,y_n))&\leq 0
\end{align*}
and choosing $\kappa>k_f$ gives a contradition.\qed\\ 

\bibliographystyle{plain}
\bibliography{dbl-obst-qvi_ref}

\begin{thebibliography}{10}

\bibitem{Bandini18}
E.~Bandini, A.~Cosso, M.~Fuhrman, and H.~Pham.
\newblock Backward sdes for optimal control of partially observed
  path-dependent stochastic systems: a control randomization approach.
\newblock {\em Ann. Appl. Probab.}, 28(3):1634--1678, 2018.

\bibitem{Barles94}
G.~Barles.
\newblock {\em Solutions de viscosit\'e des \'equations de Hamilton-Jacobi},
  volume~17 of {\em Math\'ematiques et Applications}.
\newblock Springer, Paris, 1994.

\bibitem{Barles97}
G.~Barles, R.~Buckdahn, and E.~Pardoux.
\newblock Backward stochastic differential equations and integral-partial
  differential equations.
\newblock {\em Stoch. stoch. rep.}, 60(1-2):57--83, 1997.

\bibitem{BensLionsImpulse}
A.~Bensoussan and J.L. Lions.
\newblock {\em Impulse Control and Quasivariational inequalities}.
\newblock Gauthier-Villars, Montrouge, France, 1984.

\bibitem{Bouchard09}
B.~Bouchard.
\newblock A stochastic target formulation for optimal switching problems in
  finite horizon.
\newblock {\em Stochastics}, 81:171--197, 2009.

\bibitem{Choukroun15}
S.~Choukroun, A.~Cosso, and H.~Pham.
\newblock Reflected bsdes with nonpositive jumps, and controller-and-stopper
  games.
\newblock {\em Stochastic Process. Appl.}, 125:597--633, 2015.

\bibitem{Cosso13}
A.~Cosso.
\newblock Stochastic differential games involving impulse controls and
  double-obstacle quasi-variational inequalities.
\newblock {\em SIAM J. Control Optim.}, 3(51):2102--2131, 2013.

\bibitem{UsersGuide}
M.~G. Crandall, H.~Ishii, and P.~L. Lions.
\newblock Users guide to viscosity solutions of second order partial
  differential equations.
\newblock {\em Bulletin of the American Mathematical Society}, 27(1):1--67,
  1992.

\bibitem{BollanSWG1}
B.~Djehiche, S.~Hamad\`ene, and M.~Morlais.
\newblock Viscosity solutions of systems of variational inequalities with
  interconnected bilateral obstacles.
\newblock {\em Funkcialaj Ekvacioj}, 58(1):135--175, 2015.

\bibitem{DjehicheSWG2}
B.~Djehiche, S.~Hamad\`ene, M.-A. Morlais, and X.~Zhao.
\newblock On the equality of solutions of max-min and min-max systems of
  variational inequalities with interconnected bilateral obstacles.
\newblock {\em J. Math. Anal. Appl.}, 452:148--175, 2017.

\bibitem{Dumitrescu15}
R.~Dumitrescu, M.-C. Quenez, and A.~Sulem.
\newblock Optimal stopping for dynamic risk measures with jumps and obstacle
  problems.
\newblock {\em J Optim Theory Appl}, 167:219--242, 2015.

\bibitem{ElKaroui1}
N.~El-Karoui, C.~Kapoudjian, E.~Pardoux, S.~Peng, and M.~C. Quenez.
\newblock Reflected solutions of backward {SDEs} and related obstacle problems
  for {PDEs}.
\newblock {\em Ann. Probab.}, 25(2):702--737, 1997.

\bibitem{Morlais13}
S.~Hamad\`ene and M.~A. Morlais.
\newblock Viscosity solutions of systems of pdes with interconnected obstacles
  and switching problem.
\newblock {\em Appl Math Optim.}, 67:163--196, 2013.

\bibitem{HamMor16}
S.~Hamad\`ene and M.~A. Morlais.
\newblock Viscosity solutions for second order integro-differential equations
  without monotonicity condition: the probabilistic approach.
\newblock {\em Stochastics: An International Journal of Probability and
  Stochastic Processes}, 88(4):632--649, 2016.

\bibitem{Kharroubi2010}
I.~Kharroubi, J.~Ma, H.~Pham, and J.~Zhang.
\newblock Backward sdes with constrained jumps and quasi-variational
  inequalities.
\newblock {\em Ann. Probab.}, 38(2):794--840, 2010.

\bibitem{KharroubiPham15}
I.~Kharroubi and H.~Pham.
\newblock {Feynman-Kac} representation for {Hamilton-Jacobi-Bellman} {IPDE}.
\newblock {\em Ann. Probab.}, 43(4):1823--1865, 2015.

\bibitem{OksenSulemBok}
B.~{\O}ksendal and A.~Sulem.
\newblock {\em Applied Stochastic Control of Jump Diffusions}.
\newblock Springer, 2007.

\bibitem{Stochastics23}
M.~Perninge.
\newblock Infinite horizon impulse control of stochastic functional
  differential equations driven by {L\'evy} processes.
\newblock {\em Stochastics}, pages 1--41, 2023.

\bibitem{qvi-rbsde}
M.~Perninge.
\newblock Probabilistic representation of viscosity solutions to
  quasi-variational inequalities with non-local drivers.
\newblock {\em ESAIM Control Optim.\,Calc.\,Var.}, 25:1--28, 2023.

\bibitem{PerningeJMAA23}
M.~Perninge.
\newblock Zero-sum stochastic differential games of impulse versus continuous
  control by fbsdes.
\newblock {\em J. Math. Anal. Appl.}, 527, 2023.

\bibitem{qvi-stop}
M.~Perninge.
\newblock Optimal stopping of bsdes with constrained jumps and related double
  obstacle pdes.
\newblock {\em arXiv:2402.17541}, 2024.

\bibitem{imp-stop-game}
M.~Perninge.
\newblock Optimal stopping of bsdes with constrained jumps and related zero-sum
  games.
\newblock {\em Stochastic Process. Appl.}, 173, 2024.

\bibitem{TangHou07}
S.~Tang and S.-H. Hou.
\newblock Switching games of stochastic differential systems.
\newblock {\em SIAM J. Control Optim.}, 46(3):900--929, 2007.

\end{thebibliography}
\end{document}